\documentclass[11pt,reqno,a4paper]{amsart}
\usepackage{color}
\usepackage{amssymb,amsmath}
\usepackage[latin1]{inputenc}
\usepackage[active]{srcltx}
\usepackage{graphicx}
\usepackage{exscale,relsize}
\usepackage{textgreek}
\usepackage{epsfig,graphics}
\usepackage{psfrag}
\usepackage{caption}
\usepackage{subcaption}
\usepackage[toc,page]{appendix}
\usepackage{bigints}
\def\N{\mathbb{N}}
\def\R{\mathbb{R}}
\def\m1{{I\!\!M}}

\def\ee{\`e}

\def\aa{\`a}

\renewcommand{\to}{\rightarrow}
\newcommand{\pa}{\partial}

\newcommand{\ino}{\int_{\Omega}}

\newcommand{\inpo}{\int_{\pa \Omega}}

\newcommand{\fo}{\forall}

\newcommand{\rife}[1]{(\ref{#1})}
\newcommand{\ov}[1]{\overline{#1}}

\newcommand{\scp}{\scriptstyle}
\newcommand{\sscp}{\scriptscriptstyle}
\newcommand{\dsp}{\displaystyle}

\renewcommand{\dfrac}{\displaystyle\frac}
\newcommand{\finedim}{\hspace{\fill}$\square$}
\newcommand{\intbar}{\mathop{\int\makebox(-15.5,0){\rule[6pt]{.7em}{0.3pt}}\kern-6pt}\nolimits}

\newcommand{\ii}{\infty}

\newcommand{\eps}{\varepsilon}
\newcommand{\dt}{\delta}

\newcommand{\al}{\alpha}
\newcommand{\be}{\beta}

\newcommand{\sg}{\sigma}
\newcommand{\ga}{\gamma}
\newcommand{\om}{\Omega}
\newcommand{\lm}{\lambda}

\newcommand{\rl}{\mbox{\Large \textrho}_{\!\sscp \lm}}
\newcommand{\rla}{\mbox{\Large \textrho}_{\!\sscp \lm,\sscp \al}}

\newcommand{\rlq}{(\mbox{\Large \textrho}_{\! \sscp \lm})^{\frac1q}}

\newcommand{\rlqa}{(\mbox{\Large \textrho}_{\! \sscp \lm,\al})^{\frac1q}}

\renewcommand{\rho}{\mbox{\Large \textrho}}
\newcommand{\rh}{\mbox{\Large \textrho}}

\newcommand{\pl}{\psi_{\sscp \lm}}
\newcommand{\xil}{\ul}
\newcommand{\ul}{u_{\sscp \lm}}

\newcommand{\val}{v_{\sscp I}}

\newcommand{\vxi}{\xi}
\newcommand{\wl}{w_{\sscp \lm}}

\newcommand{\vl}{\eta_{\sscp \lm}}

\newcommand{\tl}{\tau_{\ssl}}

\newcommand{\ssl}{\sscp \lm}
\newcommand{\ml}{m_{\sscp \lm}}
\newcommand{\all}{\al_{\ssl}}
\newcommand{\gal}{\ga_{\sscp I}}
\newcommand{\el}{E_{\ssl}}

\newcommand{\fbi}{{\bf (F)$_{I}$}}
\newcommand{\prl}{{\textbf{(}\mathbf P\textbf{)}_{\mathbf \lm}}}

\newtheorem{theorem}{Theorem}[section]
\newtheorem{proposition}[theorem]{Proposition}
\newtheorem{lemma}[theorem]{Lemma}
\newtheorem{corollary}[theorem]{Corollary}
\newtheorem{remark}[theorem]{Remark}
\newtheorem{definition}[theorem]{Definition}
\newcommand{\brm}{\begin{remark}\rm}
\newcommand{\erm}{\end{remark}}
\newcommand{\bdf}{\begin{definition}\rm}
\newcommand{\edf}{\end{definition}}
\newcommand{\bte}{\begin{theorem}}
\newcommand{\ete}{\end{theorem}}
\newcommand{\bpr}{\begin{proposition}}
\newcommand{\epr}{\end{proposition}}
\newcommand{\ble}{\begin{lemma}}
\newcommand{\ele}{\end{lemma}}
\newcommand{\bco}{\begin{corollary}}
\newcommand{\eco}{\end{corollary}}
\newcommand{\beq}{\begin{equation}}
\newcommand{\eeq}{\end{equation}}
\newcommand{\bdm}{\begin{displaymath}}
\newcommand{\edm}{\end{displaymath}}

\newcommand{\graf}[1]{\left\{\begin{array}{ll}#1\end{array}\right.}

\def\sideremark#1{\ifvmode\leavevmode\fi\vadjust{\vbox to0pt{\vss
 \hbox to 0pt{\hskip\hsize\hskip1em \vbox{\hsize2.1cm\tiny\raggedright\pretolerance10000 \noindent #1\hfill}\hss}\vbox to15pt{\vfil}\vss}}}

\begin{document}
\numberwithin{equation}{section}
\parindent=0pt
\hfuzz=2pt
\frenchspacing

\title[Uniqueness for free boundary problems]{
On the uniqueness and monotonicity of solutions of \\ free boundary problems}


\author[D. Bartolucci]{Daniele Bartolucci$^{(\dag)}$}
\address{Daniele Bartolucci, Department of Mathematics, University of Rome \emph{"Tor Vergata"}, Via della ricerca scientifica n.1, 00133 Roma.}
\email{bartoluc@mat.uniroma2.it}

\author[A. Jevnikar]{Aleks Jevnikar}
\address{Aleks Jevnikar, Department of Mathematics, Computer Science and Physics, University of Udine, Via delle Scienze 206, 33100 Udine, Italy.}
\email{aleks.jevnikar@uniud.it}

\thanks{$^{(\dag)}$Research partially supported by:
Beyond Borders project 2019 (sponsored by Univ. of Rome "Tor Vergata") "{\em Variational Approaches to PDE's}",
MIUR Excellence Department Project awarded to the Department of Mathematics, Univ. of Rome Tor Vergata, CUP E83C18000100006.}

\begin{abstract}
For any smooth and bounded domain $\om\subset \R^N$, we prove uniqueness of positive solutions
of free boundary problems arising in plasma physics on $\om$ in a neat interval depending only by the best constant of the
Sobolev embedding $H^{1}_0(\om)\hookrightarrow L^{2p}(\om)$, $p\in [1,\frac{N}{N-2})$ and
show that the boundary density and a suitably defined energy
share a universal monotonic behavior. At least to our knowledge, for $p>1$, this
is the first result about the uniqueness for a domain which is not a two-dimensional
ball and in particular the very first result about the monotonicity of solutions, which seems to be new even for $p=1$.
The threshold, which is sharp for $p=1$, yields a new condition
which guarantees that there is no free boundary inside $\om$.
As a corollary, in the same range, we solve a long-standing open problem (dating back to the work of Berestycki-Brezis in 1980) about the uniqueness
of variational solutions. Moreover, on a two-dimensional ball we describe the full branch of positive solutions, that is,
we prove the monotonicity along the curve of positive solutions until the boundary density vanishes.
\end{abstract}
\maketitle

\tableofcontents


\setcounter{section}{0}
\setcounter{equation}{0}
\section{\bf Introduction}

Letting $\om\subset \R^N$, $N\geq2$, be an open and bounded domain of class $C^{3}$, we are concerned with the free boundary problem
$$
\graf{-\Delta v = (v)_{+}^p\quad \mbox{in}\;\;\om\\ \\
-\bigintss\limits_{\pa\om} \dfrac{\pa v}{\pa\nu}=I \\ \\
v=\ga \quad \mbox{on}\;\;\pa\om
}\qquad \qquad \mbox{\bf (F)}_{I}
$$
for the unknowns $\ga \in \R$ and
$v\in C^{2,r}(\ov{\om}\,)$, $r\in (0,1)$. Here $(v)_+$ is the positive part of $v$, $\nu$ is the exterior unit normal, $I> 0$ and $p\in (1,p_{N})$ are fixed, with
$$
p_{ N}=\graf{+\ii,\;N=2 \vspace{0.2cm}\\ \frac{N}{N-2}\,,\; N\geq 3.}
$$
Up to a suitable scaling, one can
assume without loss of generality that $|\om|=1$ and peak any positive constant to multiply $(v)_+^p$.\\

The study of \fbi\, is motivated by Tokamak's plasma physics and we refer to
\cite{FL,Te} and to \cite{Kad,Mer} for a physical description of the problem.
Initiated in \cite{Te,Te2} with the model case $p=1$ (which requires a slightly different formulation, see Remark \ref{remeq}), a systematic analysis of \fbi\, has been undertaken in \cite{BeBr}
with more general operators and
the nonlinearity $(v)_{+}^p$ being replaced by
the positive part of any continuous and increasing function with growth
of order at most $p$, for some $p\in [1,p_{N})$. Indeed, the threshold $p_{N}$ turns out to be a natural critical exponent for \fbi, see \cite{BeBr,Ort} and the discussion in the sequel.\\
However, among many other things, it has been shown in \cite{BeBr} that for any $I>0$ there exists at least one solution
of \fbi. A lot of work has been done to understand solutions of \fbi\, for
$p\in (1,p_{N}]$, see \cite{AmbM,BMar,BMar1,BSp,FW,Kor3,LP,Liu,Mar,Ort2,suzt,W,WY,We,Wo}, and in the model case $p=1$, see
\cite{CF,CPY,dam,Fr,Gal,Pu1,Pudam,Sch,Shi}, and the references quoted therein. Although we will not discuss this point here, a lot of work has been done
in particular to understand the regularity (for solutions with $\ga<0$) of the free boundary $\pa\{x\in\om\,|\,v(x)>0\}$,
see for example \cite{FL,JP,KNS,KS}. Finally, for a closely related problem see also \cite{BGG,CYYZ,SV-S} and references quoted therein.\\

However, the uniqueness results at hand about \fbi\, seem to concern so far either
the model case $p=1$ and $N\geq 2$ (\cite{Gal,Pudam,Te2}) or the case of the ball for $N=2$ and $p>1$ (\cite{BSp}), or either
the case $N\geq 2$ where $(v)^p_+$ is replaced by the positive part of any continuous and increasing nonlinearity
satisfying a suitable uniformly Lipschitz condition (whence at most of linear growth), see \cite{BeBr}. Otherwise,
we are not aware of any uniqueness result for solutions of \fbi\,
neither for positive solutions, nor for the so-called variational solutions, see \cite{BeBr} and \cite{BMar}.
Actually, in dimension $N\geq 3$ and $p=\frac{N+2}{N-2}$
uniqueness fails for \fbi\, on a ball for $I$ small enough, see \cite{BMar1}. On the other hand, a non-uniqueness result for \fbi, with the nonlinearity rewrited in the form $\lm(v)_{+}^p$, is obtained in \cite{CF,CPY,Liu,Sch} for any $p\in[1,p_N)$ and $\lm$ sufficiently large.
However, our main motivation comes from the fact that, at least to our knowledge, no results at all are available
so far about the shape of branches of solutions $(\gal,\val)$ of \fbi.\\

The main idea here is to identify the
natural spectral setting together with the quantities which should share nice monotonicity properties (see also \cite{B2,BJ,BW} where a different class of problems is considered). The first step is to move to a dual formulation (\cite{BeBr,Te2}) of \fbi\,
via the constrained problem,
$$
\graf{-\Delta \psi =(\al+{\lm}\psi)^p\quad \mbox{in}\;\;\om\\ \\
\bigintss\limits_{\om} {\dsp \left(\al+{\lm}\psi\right)^p}=1\\ \\
\psi>0 \quad \mbox{in}\;\;\om, \quad \psi=0 \quad \mbox{on}\;\;\pa\om \\ \\
\al\geq0
}\qquad \prl
$$
for the unknowns $\al\in\R$ and $\psi \in C^{2,r}_{0,+}(\ov{\om}\,)$. Here, $\lm\geq 0$ and ${p\in [1,p_N)}$ are fixed and for $r\in (0,1)$ we set
$$
C^{2,r}_0(\ov{\om}\,)=\{\psi \in C^{2,r}(\ov{\om}\,)\,:\, \psi=0\mbox{ on }\pa \om\},\;
C^{2,r}_{0,+}(\ov{\om}\,)=\{\psi \in C^{2,r}_0(\ov{\om}\,)\,:\, \psi> 0\mbox{ in } \om\}.
$$

It is useful to define positive solutions as follows.\\

{\bf Definition.}
{\it We say that $(\gal,\val)$ is a positive solution
of \fbi\, if $\gal> 0$.}
{\it We say that $(\all,\pl)$ is a positive solution of {\rm $\prl$} if $\all>0$.}

\bigskip

We remark that since $|\om|=1$ and $\lm\geq 0$ by assumption, then if $(\all,\pl)$ solves $\prl$ we necessarily have,
$$
\all\leq 1,
$$
and the equality holds if and only if $\lm=0$. We will frequently use this fact without further comments.
Actually, if $\lm=0$, then $\prl$ admits a unique solution $(\al_{\sscp 0}, \psi_{\sscp 0})=(1,G[1])$ satisfying
$$
\graf{-\Delta \psi_{\sscp 0} =1\quad \mbox{in}\;\;\om\\ \\
\psi_{\sscp 0}=0 \quad \mbox{on}\;\;\pa\om
}
$$
where we define,
$$
G[\rho](x)=\ino G_{\om}(x,y)\rho(y)\,dy,\;x\in\om.
$$
Here $G_{\om}$ is the Green function of $-\Delta$ with Dirichlet boundary conditions on $\om$.
Obviously, to say that $(\all,\pl)$ is a solution of $\prl$ is the same as to say that $\pl=G[\rl]$ and $\ino \rl=1$,
where, here and in the rest of this work, we set
$$
\rl=(\all+\lm\pl)^p.
$$

\brm\label{remeq} {\it
Let $q$ be the conjugate index of $p$, that is
$$
\frac1p+\frac1q=1.
$$
For any fixed $\lm>0$ and $p>1$, $(\all,\pl)$ is a positive
solution of {\rm $\prl$} if and only if,
for $I=I_{\ssl}=\lm^{q}$, $(\gal,\val)=(\lm^{\frac{1}{p-1}}\all,\lm^{\frac{1}{p-1}}(\all+\lm\pl))$ is a
positive solution of \fbi. Therefore,
in particular, if $(\gal,\val)$ solves \fbi\, then
$(\all,\pl)=(I^{-\frac{1}{p}}\gal,I^{-1}(\val -\gal))$ solves
{\rm $\prl$} and the identity $I^{-\frac{1}{p}}\val=\all+\lm\pl$ holds. This correspondence is singular
for $p=1$ and in this case we will stick with the formulation {\rm $\prl$}. Indeed, it is readily seen that for $p=1$, {\rm $\prl$} is equivalent to a more general problem than \fbi\, and positive solutions of {\rm $\prl$} correspond to positive solutions of \fbi\, where the first equation is replaced by $-\Delta v = \lm(v)_{+}$.}
\erm

\medskip

All the statements below are concerned with $\prl$, the corresponding statements
about \fbi\, being immediately recovered via the latter remark. For the sake of clarity, we point out that if a map $\mathcal{M}$ from an interval $[a,b]\subset \R$ to a Banach space $X$ is said to be real analytic,  then it is understood that $\mathcal{M}$ can be extended in an open neighborhood of $a$ and $b$ where it admits a power series expansion, totally convergent in the $X$-norm.\\

For fixed $t\geq 1$ and $|\om|=1$ we denote
$$
\Lambda(\om,t)=\inf\limits_{w\in H^1_0(\om), w\equiv \!\!\!\!/ \;0}\dfrac{\ino |\nabla w|^2}{\left(\ino |w|^{t}\right)^{\frac2t}}\,,
$$
which provides the best constant in the Sobolev embedding $\|w\|_p\leq S_p(\om)\|\nabla w\|_2$, $S_p(\om)=\Lambda^{-2}(\om,p)$,
$p\in[1,2p_{N})$.
For $(\all,\pl)$ a solution of $\prl$ we define the energy,
$$
\el:=\frac12 \ino\rl \pl\equiv \frac12 \ino |\nabla \pl|^2,
$$
see Appendix \ref{appD} for more comments about the latter definition.\\

Set
$$
\lm^*(\om,p)=\sup\{\lm>0\,:\,\al_{\mu}>0\,\mbox{\rm \,for any solution of }
{\textbf{(}\mathbf P\textbf{)}_{\mathbf \mu}},\,\forall\,\mu<\lm\}.
$$
It is not difficult to show (see Lemma \ref{lmsmall} in Appendix \ref{appF}) that the latter quantity is well defined. Finally, we denote by $\mathcal{G}(\om)$ the set of solutions of $\prl$ for $\lm\in [0,\frac1p\Lambda(\om,2p))$, $p\in[1,p_{N})$ and
let $B_r=\{x\in \R^N\,:\, |x|<r\}$ with volume $|B_r|$ and
$\mathbb{D}_{\sscp N}$ be the $N$-dimensional
ball of unit volume.
\bte\label{thmLE}
Let $p\in [1,p_{N})$. Then $\lm^*(\om,p)\geq \frac1p \Lambda(\om,2p)$ and the equality holds if and only if $p=1$. Moreover, we have:
\begin{itemize}
\item[1.] \emph{(Uniqueness):} for any $\lm \in [0,\frac1p \Lambda(\om,2p))$
there exists a unique solution $(\all,\pl)$ of {\rm $\prl$}.

\item[2.] \emph{(Monotonicity):} $\mathcal{G}(\om)$ is a real analytic simple curve of positive solutions $[0,\frac1p\Lambda(\om,2p))\ni \lm\mapsto (\all,\pl)$ such that, for any $\lm\in [0,{\textstyle \frac1p}\Lambda(\om,2p))$,
$$
\frac{d \all}{d\lm}<0\quad\mbox{ and }\quad \frac{d \el}{d\lm}>0,
$$
and
$$
\all=1+\mbox{\rm O}(\lm),\;\;\pl=\psi_{\sscp 0}+\mbox{\rm O}(\lm),
\;\; \el=E_0(\om)+\mbox{\rm O}(\lm), \quad \mbox{ as }\lm\to 0^+,
$$
where,
$$
E_{0}(\om)=\frac12\ino \ino G_{\om}(x,y)\,dxdy\leq
E_{0}(\mathbb{D}_{\sscp N})=
{\textstyle\dfrac{|B_1|^{\scp -\frac 2 N}}{4(N+2)}}.
$$
In particular $\mathcal{G}(\om)$ can be extended continuously on $\lm \in [0,\frac1p \Lambda(\om,2p)]$ with
$(\all,\pl)\to (\ov{\al},\ov{\psi})$ as $\lm\to \frac1p\Lambda(\om,2p)^-$ and $\ov{\al}=0$ if and only if $p=1$.
\end{itemize}
\ete

\medskip

The threshold
${\textstyle \frac1p}\Lambda(\om,2p)$ is sharp for $p=1$ since in this case it is well known (\cite{BeBr,Pudam,Te2}) that solutions
are unique and positive if and only if
$\lm<\lm^*(\om,1)=\lm^{(1)}(\om)$ while clearly $\Lambda(\om,2)=\lm^{(1)}(\om)$, where
$\lm^{(1)}(\om)$ is the first eigenvalue of $-\Delta$ on $\om$ with Dirichlet boundary conditions. The sharpness of the above result, for $p>1$, will be discussed in a subsequent work. The bound $\lm^*(\om,p)\geq \frac1p\Lambda(\om,2p)$ is a superlinear generalization of
the sufficient conditions, obtained in the sublinear case in \cite{AmbM,BMar1},
which guarantee that for variational solutions there is no free boundary inside $\om$, see Theorem \ref{thmvar}. Clearly, $2E_0(\om)$ is just the torsional rigidity of $\om$.\\
At least to our knowledge, for $p>1$, this is the first result about the uniqueness of solutions for a domain which is not a two-dimensional
ball. In particular, this is the very first result about the qualitative behavior of the branch of solutions and seems to be new even for $p=1$, in which case, we actually describe the full branch of positive solutions, i.e. until the boundary density vanishes. Both the uniqueness and monotonicity of solutions hold for any smooth and bounded domain, in any dimension and for any subcritical exponent.\\

Therefore, we succeed in the construction of
a branch of positive solutions, that is, the case where there is no free boundary inside $\om$,
emanating from the unique positive solution $(\al_{\sscp 0},\psi_{\sscp 0})$
which share a universal (i.e. independent of $\om$) monotonic profile,
and prove that they are also the unique solutions of $\prl$ in $[0,{\textstyle \frac1p}\Lambda(\om,2p))$.
The corresponding branch of solutions for \fbi\,, $p\in (1,p_{N})$, is real analytic in $(0,(\frac1p\Lambda(\om,2p))^{q})$,
continuous in $(0,(\frac1p\Lambda(\om,2p))^{q}]$ and tends to the trivial solution $(\ga_0,v_0)\equiv (0,0)$ as $I\to0^+$.
We remark that, while $\all=I^{-\frac1p}\gal$ is decreasing,
$\gal$ need not be monotone and indeed it is not, at least on $\om=\mathbb{D}_{\sscp 2}$, see Remark \ref{remdisk} below.\\

Theorem \ref{thmLE} is also relevant in the study of the so-called
variational solutions of \fbi, see \cite{BeBr} and \cite{BMar,BMar1}.
A solution $(\gal,\val)$ of \fbi\, is a variational solution of \fbi\, if it is a minimizer of
the functional,
\beq\label{psiI}
\Psi_{I}(v)=\frac12\ino |\nabla v|^2-\frac{1}{p+1}\left(\ino (v)_+^{p+1}\right)+ Iv(\pa \om),
\eeq
on
$$
\mathcal{H}_I=\left\{v\in H\,|\,\ino (v)_+^p=I\right\},
$$
where $I>0$ and $H$ is the space of $H^{1}(\om)$ functions whose boundary trace is constant. Analogously, variational solutions of $\prl$ are related to solutions of the dual variational principle {\bf (VP)}, see Appendix \ref{appD} for more details.
It has been shown in \cite{Te2} that for $p=1$ and any $\lm>0$ there exists at least one variational solution of $\prl$, while in \cite{BeBr} the authors proved that for any $p\in (1,p_N)$ and any $I>0$ there exists at least one variational solution of \fbi.
Actually, we know from \cite{BMar1,Gal,Pudam,Te2} that the following holds true.

\medskip

{\bf Theorem A.} (\cite{BMar1,Gal,Pudam,Te2})\\
1. {\it Let $p\in(1,p_{N})$ and $(\gal,\val)$ be a variational solution of \fbi.
Then there exists $I^{**}(\om,p)\in (0,+\ii)$ such that $\gal>0$ if and only if $I\in (0,I^{**}(\om,p))$.}

\medskip

2. {\it Let $p=1$ and $(\all,\pl)$ be a any solution of \emph{$\prl$}. Then there exists $\lm^{**}(\om,1)\in (0,+\ii)$ such that $\all>0$ if and only if $\lm\in (0,\lm^{**}(\om,1))$. Moreover, $\lm^{**}(\om,1)=\lm^{(1)}(\om)$.}\\

For $p>1$ there exists a one to one correspondence (\cite{BeBr}) between variational solutions of
\fbi\, and variational solutions of $\prl$ and in particular we define,
$$
\lm^{**}(\om,p)=(I^{**}(\om,p))^{\frac1q}, \quad p>1,
$$
which shares the same property of $I^{**}(\om,p)$ in Theorem A, see Corollary \ref{eqlms} in Appendix \ref{appD}.
By the uniqueness results for the model case $p=1$ (\cite{Gal,Pudam,Te2}), any positive solution is a variational solution and then in particular we have $\lm^{**}(\om,1)=\lm^{*}(\om,1)$. Similarly, the uniqueness property in \cite{BSp} for $\om=\mathbb{D}_{\sscp 2}$ yields that
$\lm^{**}(\mathbb{D}_{\sscp 2},p)=\lm^{*}(\mathbb{D}_{\sscp 2},p)$.
Excluding these model cases, uniqueness was not known so far neither for variational solutions. Concerning this point
we have the following,
\bte\label{thmvar}
For $p\in[1,p_{N})$ and for any $\lm\in [0,{\textstyle \frac1p}\Lambda(\om,2p))$ there exists a unique variational solution
$(\all,\pl)$ of {\rm $\prl$}. For $p\in(1,p_{N})$ and for any $I\in (0,({\textstyle \frac1p}\Lambda(\om,2p))^{q})$
there exists a unique variational solution
$(\gal,\val)$ of \fbi. In particular the set of
unique variational solutions of {\rm $\prl$} in $[0,{\textstyle \frac1p}\Lambda(\om,2p))$ coincides with $\mathcal{G}(\om)$ and the following
inequalities hold,
\beq\label{istest}
\lm^{**}(\om,p)\geq \lm^*(\om,p)\geq \frac1p\Lambda(\om,2p)\geq \frac1p\Lambda(\mathbb{D}_{\sscp N},2p).
\eeq
\ete

\medskip

Observe that from Theorem \ref{thmLE} and
\rife{istest} we have $\lm^{**}(\om,p)\geq \lm^*(\om,p)>\frac1p\Lambda(\om,2p)$, for $p>1$. On the other hand, the set of positive variational solutions is
not empty for any $\lm\in (\frac1p\Lambda(\om,2p), \lm^{**}(\om,p))$, see \cite{BMar,BeBr}. The continuation of the curve of solutions $\mathcal{G}(\om)$, under generic assumptions, beyond $\frac1p\Lambda(\om,2p)$, enjoying uniqueness and monotonicity properties, will be the topic of a future paper.\\

Concerning this point, for $\om=\mathbb{D}_{\sscp 2}$ and for any $p\in[1,+\infty)$, we succeed in the description of the shape and monotonicity
of the branch of positive solutions with $\lm<{\lm}^*(\mathbb{D}_{\sscp 2},p)$. We point out that this is the full branch of positive solutions since, as already observed, $\all>0$ if and only if $\lm<{\lm}^*(\mathbb{D}_{\sscp 2},p)$. The sharp positivity threshold has been recently explicitly computed, see Theorem 1.2 in \cite{BJ2}, and takes the form
$$
{\lm}^*(\mathbb{D}_{\sscp 2},p)=
\left(\frac{8\pi}{p+1}\right)^{\frac{p-1}{2p}}\Lambda^{\frac{p+1}{2p}}(\mathbb{D}_{\sscp 2},p+1).
$$ 
We point out that here only uniqueness was known so far (\cite{BSp}). Let $\mathcal{G}^{*}(\mathbb{D}_{\sscp 2})$ be the set of unique
solutions of $\prl$, then we have,
\bte\label{tdisk}
Let $p\in[1,+\infty)$. Then, $\mathcal{G}^{*}(\mathbb{D}_{\sscp 2})$ is a continuous curve, defined in $[0, {\lm}^*(\mathbb{D}_{\sscp 2},p)]$, such that,
$$
\all\mbox{ is strictly decreasing and } \el\mbox{ is strictly increasing in }(0,{\lm}^*(\mathbb{D}_{\sscp 2},p)),
$$
and
$$
\al_{\sscp \lm^*(\mathbb{D}_{\sscp 2},p)}= 0,\quad E_{\lm^*(\mathbb{D}_{\sscp 2},p)}=\frac{p+1}{16\pi}.
$$
\ete

\medskip

\brm\label{remdisk}{\it Actually, we prove something more about the regularity of $\all, \el$, see section \ref{sec6}. Observe that, along $\mathcal{G}^{*}(\mathbb{D}_{\sscp 2})$, by the monotonicity of $\el$, we have $\el\geq E_0(\mathbb{D}_{\sscp 2})=\frac{1}{16\pi}$. On the other hand, by Theorem 1.1 in \cite{BJ2}, we know $\el\leq\frac{p+1}{16\pi}$. Thus, along the branch, the solutions span the full energy range $[\frac{1}{16\pi},\frac{p+1}{16\pi}]$.
Besides, for $p>1$, the parameter $\gal=\lm^{\frac{1}{p-1}}\all$ varies continuously from
$\gamma_{0}=0$ to $\gamma_{I^*(\mathbb{D}_{\sscp 2},p)}=0$, whence in particular it is not monotone
along the branch. For $\gal$ small enough there are at least two distinct pairs $\{I,\val\}$ solving \fbi.}
\erm

\medskip

The proof of Theorem \ref{thmLE} is based on the
bifurcation analysis of the vectorial solutions $(\all,\pl)$ of the constrained problem $\prl$.\\
The crucial property of $\prl$ relies on its linearized operator and on a suitably defined constrained spectral theory, see the definition \rife{eLl}
of $L_{\ssl}$ and the related eigenvalues equation \rife{lineq0} in section \ref{sec2}.
The use of this operator is rather delicate since $L_{\ssl}$ arises as the linearization of a constrained problem
$(\ino \rl=1)$ with respect to $(\all,\pl)$, which yields a non-local problem. As a consequence
it is not true in general that its first eigenvalue, which we denote by $\sg_1(\all,\pl)$, see \rife{4.1},
is simple and neither that
if $\sg_1(\all,\pl)$ is positive then the maximum principle holds. For example this is exactly
what happens for $\lm=0$ on $\mathbb{D}_{\sscp 2}$, where $\sg_1(\al_{\sscp 0},\psi_{\sscp 0})$ can be
evaluated explicitly (see \cite{BJ}) and one finds that
$\sg_1(\al_{\sscp 0},\psi_{\sscp 0})=\lm^{(2,0)}(\mathbb{D}_{\sscp 2})\simeq \pi (3,83)^2$ has three eigenfunctions,
two of which indeed change sign. Here
$\lm^{(2,0)}(\om)>\lm^{(1)}(\om)$ is the first non vanishing eigenvalue of $-\Delta$ on $\om$ on
the space of $H^1(\om)$ vanishing mean functions
with constant boundary trace.\\
In few words we argue as follows. We first prove that,
if for a positive solution $(\all,\pl)$ with $\lm\geq 0$ we have $0\notin \sigma(L_{\ssl})$, where $\sigma(L_{\ssl})$ stands for the spectrum of $L_{\ssl}$, then the set of solutions of $\prl$
is locally a real analytic curve of positive solutions. In particular a real analytic curve of
positive solutions exists around $(\al_{\sscp 0},\psi_{\sscp 0})$.
A crucial information is that if $\sg_1(\all,\pl)$ is positive,
then $\frac{d \all}{d\lm}<0$ and $\frac{d \el}{d\lm}>0$, see section \ref{sec3}.
We evaluate the sign of the derivative of $\el$ by a careful analysis of the Fourier modes of
$L_{\ssl}$. On the other hand, the monotonicity of $\all$ is hidden in tricky identities involving the solution $(\all,\pl)$. At this point one can see that if for some $\ov{\lm}>0$ it holds $\sg_1(\all,\pl)>0$ for any solution
$(\all,\pl)$ of $\prl$ with $\lm\leq \ov{\lm}$ and if $\ov{\al}>0$ for a solution at $\ov{\lm}$,
then there exists one and only one solution of $\prl$ in $[0,\ov{\lm}\,]$ which, by the monotonicity of $\all$,
is also a positive solution, see Lemma \ref{thm05}.
This interesting property seems to be the generalization to free boundary problems of similar properties
for minimal solutions of semilinear elliptic equations with strictly positive, increasing and convex nonlinearities
(\cite{CrRab,KK}).
In particular, as a consequence of these facts, the proof of Theorem \ref{thmLE}
is reduced to an a priori bound from below away from zero for $\sg_1(\all,\pl)$ and $\all$.\\
The spectral estimate is obtained by using that $\sg_1(\all,\pl)$ is strictly larger
than the "standard" eigenvalue $\nu_{1,\lm}$ corresponding to the problem without constraints, see \rife{xi1s}.
This is crucial as it rules out at once the natural strong competition between the value of the $\lm$-threshold and the
bound from below for $\all$.
On the other side, the fact that $\lm^*(\om,p)\geq \frac1p\Lambda (\om,p)$ follows essentially
as consequence of the fact that if $(0,\pl)$ solves $\prl$, then $\lm\pl$ is proportional to a
positive strict subsolution of the corresponding linearized problem.\\

Concerning Theorem \ref{tdisk} we remark that
even on $\mathbb{D}_{\sscp 2}$ it is not trivial to catch the shape of the full branch of positive solutions.
Indeed, in principle one should handle the presence of vanishing higher order eigenvalues along the branch of solutions $\mathcal{G}(\mathbb{D}_{\sscp 2})$. Exploiting the uniqueness of solutions (\cite{BSp}) together with a characterization of variational solutions, see Lemma \ref{prmin}, we are led to study the eventuality of singularities at $\sg_1(\all,\pl)=0$. The non-local structure of $L_{\ssl}$
makes the bifurcation analysis rather delicate. We handle this problem by a generalization of the
Crandall-Rabinowitz (\cite{CrRab}) bending result suitable to describe solutions $(\all,\pl)$ of $\prl$ on
any domain $\om$ near a singular point satisfying a suitable transversality condition.
This result, although we will not state it explicitly, is of independent interest as it
is just, via Remark \ref{remeq}, a Crandall-Rabinowitz bending-type result suitable to be applied to
positive solutions of \fbi. We stress that the latter result would not work out if we would stick with the standard spectral analysis as it exploits the modified (constrained) spectral setting in its full strength.
Besides its application to the case $\om=\mathbb{D}_2$, it will be exploited in a future work
to obtain sufficient conditions to continue $\mathcal{G}(\om)$,
under generic assumptions, until the boundary density vanishes.\\

The paper is organized as follows. In section \ref{sec2} we set up the spectral and bifurcation analysis. Then, in section \ref{sec3} we show the monotonicity of solutions while in section \ref{sec4} we prove the uniqueness result. Finally, in section \ref{sec6} we discuss the two-dimensional ball case. A brief discussion about variational solutions and uniqueness of solutions for $\lm$ small is postponed to Appendixes \ref{appD} and \ref{appF}, respectively.

\bigskip
\bigskip

\section{\bf Spectral and bifurcation analysis}\label{sec2}
In this section we develop the spectral and bifurcation analysis for solutions of $\prl$ with $\lm\geq 0$ and
$p\in[1,p_{N})$. Lemma \ref{lem1.1} is the basic result which we will use to describe branches of solutions
of $\prl$ at regular points. We will also prove Lemma \ref{prmin} and Theorem \ref{spectral1} which will be needed in section
\ref{sec6} for the two-dimensional ball case. In particular we will need there
a generalization of the bending result of \cite{CrRab} for solutions of $\prl$, see Proposition \ref{pr3.1} below.\\

From now on,  we will denote
$$
\tl=p\lm.
$$
Also, we write $\rla=\rla(\psi)=(\al+\lm\psi)^p$ and $\rl=(\all+\lm\pl)^p$, and then in particular,
$$
(\rla(\psi))^{\frac1q}=(\al+\lm\psi)^{p-1} \mbox{ and } \rlq=(\all+\lm\pl)^{p-1},
$$
where it is understood that $(\rla(\psi))^{\frac1q}\equiv 1 \equiv \rlq$ for $p=1$.
Whenever $(\all,\pl)$ is a solution of $\prl$ we denote,
$$
<\eta>_{\ssl}=\dfrac{\ino \rlq  \eta}{\ino \rlq}\quad \mbox{and }\quad  [\eta]_{\ssl}=\eta \,-<\eta>_{\ssl},
$$
and define,
$$
<\eta,\varphi>_{\ssl}:=\dfrac{\ino \rlq  \eta\varphi}{\ino \rlq}\quad \mbox{and}\quad
\|\varphi\|_{\ssl}^2:=<\varphi,\varphi>_{\ssl}=<\varphi^2>_{\ssl}=\dfrac{\ino \rlq  \varphi^2}{\ino \rlq}\,,
$$
where $\{\eta,\varphi\}\subset L^2(\om)$. Clearly, for any solution $\rl$ is strictly positive in ${\om}$.
Therefore, it is easy to see that $<\cdot,\cdot>_{\ssl}$ defines a scalar product on $L^2(\om)$ whose norm is  $\|\cdot\|_{\ssl}$.
We will also adopt sometimes the useful shorthand notation,
$$
\ml=\ino \rlq.
$$

\brm\label{rem1}{\it
We will use the fact that,{\rm
$$
{\ino \rlq  [\eta]_{\ssl}^2}=\ino \rlq  (\eta\,-<\eta>_{\ssl})^2\geq 0,
$$
}where the equality holds if and only if $\eta$ is constant, whence in particular, in case $\eta\in H^{1}_0(\om)$,
if and only if $\eta$ vanishes identically. Also, since obviously $<[\eta]_{\ssl}>_{\ssl}=0$, then,
$$
<[\varphi]_{\ssl}, \eta>_{\ssl}=<\varphi, [\eta]_{\ssl}>_{\ssl}= <[\varphi]_{\ssl}[\eta]_{\ssl}>_{\ssl},\;\fo\,\{\eta,\varphi\}\subset L^2(\om).
$$
We will often use these properties when needed without further comments.}
\erm

\bigskip
In the sequel we aim to describe possible branches of solutions of $\prl$ around a positive solution, i.e. with $\all>0$. To this end, it is not difficult to construct an open subset $A_{\sscp \om}$ of the Banach space of triples $(\lm,\al,\psi)\in \R\times \R\times  C^{2,r}_0(\ov{\om}\,)$ such that, on $A_\om$, the density $\rla=\rla(\psi)=(\al+\lm\psi)^p$ is well defined and
\begin{equation} \label{bound-a}
\all+\lm\pl\geq \frac{\all}{2} \quad \mbox{ in }\ov{\om}
\end{equation}
in a sufficiently small open neighborhood in $A_{\sscp \om}$ of any triple of the form $(\lm,\all,\pl)$ whenever $(\all,\pl)$ is a positive solution of {\rm $\prl$}.\\
At this point we introduce the map,
\beq\label{eF}
F: A_{\sscp \om} \to  C^{r}(\ov{\om}\,),\quad F(\lm,\al,\psi):=-\Delta \psi -\rla(\psi).
\eeq
and
\beq\label{ePhi}
\Phi:A_{\sscp \om} \to  \R\times C^{r}(\ov{\om}\,),\quad
\Phi(\lm,\al,\psi)=\left(\begin{array}{cl}F(\lm,\al,\psi)\\ \\-1+\ino \rla\end{array}\right),
\eeq
and, for a fixed $(\lm,\al,\psi)\in A_{\sscp \om}$,
its differential with respect to $(\al,\psi)$, that is the linear operator,
$$
D_{\al,\psi}\Phi(\lm,\al,\psi):\R \times C^{2,r}_0(\ov{\om}\,) \to  \R\times C^{r}(\ov{\om}\,),
$$
which acts as follows,
$$
D_{\al,\psi}\Phi(\lm,\al,\psi)[s,\phi]=
\left(\begin{array}{cr}
D_\psi F(\lm,\al,\psi)[\phi]+d_{\al}F(\lm,\al,\psi) [s]\\ \\
\ino \left(D_\psi \rla[\phi]+d_{\al}\rla [s]\right)\end{array}\right),
$$
where we have introduced the differential operators,
\beq\label{lin}
D_\psi F(\lm,\al,\psi)[\phi]=
-\Delta \phi -\tl \rlqa \phi, \quad \phi \in C^{2,r}_0(\ov{\om}\,),
\eeq

\beq\label{lin1}
D_\psi \rla[\phi]=\tl\rlqa \phi, \quad \phi \in C^{2,r}_0(\ov{\om}\,),
\eeq
and

\beq\label{lin2}
d_{\al}F(\lm,\al,\psi) [s]=-p \rlqa s, \quad s \in \R,
\eeq

\beq\label{lin3}
d_{\al}\rla [s]= p\rlqa s, \quad s \in \R,
\eeq
where we recall $\tl=p\lm$.\\
By the construction of $A_{\sscp \om}$, see in particular \eqref{bound-a}, relying on
known techniques about real analytic functions on Banach spaces {\rm (}\cite{but}{\rm )},
it is not difficult to show that $\Phi(\lm,\al,\psi)$ is jointly real analytic in an open neighborhood of $A_{\sscp \om}$
around any triple of the form $(\lm,\all,\pl)$ whenever $(\all,\pl)$ is a positive solution of {\rm $\prl$}.\\

For fixed $\lm\geq 0$ and $p\in[1,p_{N})$, the pair $(\all,\pl)$ solves $\prl$ in the classical sense as defined in the introduction
if and only if
$\Phi(\lm,\all,\pl)=(0,0)$,  and we define the linear operator,
\beq\label{eLl}
L_{\ssl}[\phi]=-\Delta \phi-\tl \rlq [\phi]_{\ssl}.
\eeq

We say that $\sg=\sg(\all,\pl)\in\R$ is an eigenvalue of $L_{\ssl}$ if the equation,
\beq\label{lineq0}
-\Delta \phi-\tl \rlq [\phi]_{\ssl}=\sg\rlq [\phi]_{\ssl},
\eeq
admits a non-trivial weak solution $\phi\in H^1_0(\om)$. Let us define the Hilbert space,
\beq\label{y0}
Y_0:=\left\{ \varphi \in \{L^2(\om),<\cdot,\cdot >_{\ssl}\}\,:\,<\varphi>_{\ssl}=0\right\},
\eeq
and $T(f):=G[ \rlq f]$, for $f\in L^2(\om)$.
Since $T(Y_0)\subset W^{2,2}(\om)$, then the linear operator,
\beq\label{T0}
T_0:Y_0\to Y_0,\;T_0(\varphi)=G[\tl \rlq\varphi]-<G[\tl \rlq\varphi]>_{\ssl},
\eeq
is compact. By a straightforward evaluation we see that $T$ is also self-adjoint. As a consequence,
standard results concerning the spectral decomposition of self-adjoint, compact, linear operators on
Hilbert spaces show that $Y_0$ is the Hilbertian direct sum of the eigenfunctions of $T_0$, which can be represented as
$\varphi_k=[\phi_k]_{\ssl}$, $k\in\N=\{1,2,\cdots\}$,
$$
Y_0=\overline{\mbox{Span}\left\{[\phi_k]_{\ssl},\;k\in \N\right\}},
$$
for some $\phi_k\in H^1_0(\om)$, $k\in\N=\{1,2,\cdots\}$. In fact, the eigenfunction $\varphi_k$,
whose eigenvalue is $\mu_k\in \R\setminus\{0\}$, satisfies,
$$
\mu_k\varphi_k=\left(G[\tl \rlq\varphi_k]-< G[\tl \rlq\varphi_k]>_{\ssl}\right).
$$
In other words, by defining,
$$
\phi_k:=(\tl+\sg_k) G[\rlq\varphi_k],
$$
it is easy to see that $\varphi_k$ is an eigenfunction of $T_0$ with eigenvalue $\mu_k=\frac{\tl}{\tl+\sg_k}\in \R\setminus\{0\}$ if and only if
$\phi_k\in H^1_0(\om)$ and weakly solves,
\beq\label{lineqm}
-\Delta \phi_k= (\tl+\sg_k)  \rlq [\phi_k]_{\ssl}\quad\mbox{ in }\quad \om.
\eeq

In particular we will use the fact that $\varphi_k=[\phi_k]_{\ssl}$ and
\beq\label{lineq9m}
\phi_k=(\tl+\sg_k) G[\rlq[\phi_k]_{\ssl}],\quad k\in\N=\{1,2,\cdots\}.
\eeq

At this point, standard arguments in the calculus of variations show that,
\beq\label{4.1}
\sg_1=\sg_1(\all,\pl)=\inf\limits_{\phi \in H^1_0(\om)\setminus \{0\}}
\dfrac{\ino |\nabla \phi|^2 - \tl \ino \rlq [\phi]_{\ssl}^2 }{\ino \rlq [\phi]_{\ssl}^2}.
\eeq
The ratio in the right hand side of \rife{4.1} is well defined because of
Remark \ref{rem1} and in particular we have,
\beq\label{lineq2}
\tl +\sg_1>0.
\eeq
By the Fredholm alternative,
if $0\notin \{\sg_j\}_{j\in \N}$, then $I-T_0$ is an isomorphism of $Y_0$ onto itself. Clearly, we can construct an orthonormal base
of  eigenfunctions satisfying,
\beq\label{orto}
<[\phi_i]_{\ssl}, [\phi_j]_{\ssl}>_{\ssl}=0,\quad \forall i\neq j.
\eeq
The following easy to verify inequality about the standard first eigenvalue $\nu_{1,\lm}$ holds,
$$
\nu_{1,\lm}=\inf\limits_{w\in H^{1}_0(\om)\setminus\{0\}}\dfrac{\ino |\nabla w|^2-\tl \ino \rlq w^2  }
{\ino \rlq w^2}\leq \sg_1(\all,\pl).
$$

Actually, we will need the following slightly refined inequality which holds for any $p\in [1,p_{N})$,

\beq\label{xi1s}
\sg_1(\all,\pl)>\nu_{1,\lm}.
\eeq
Indeed we have,
$$
\dfrac{\ino |\nabla \phi|^2 - \tl \ino \rlq [\phi]_{\ssl}^2 }{\ino \rlq [\phi]_{\ssl}^2}=
\dfrac{\ino |\nabla \phi|^2 }{\ino \rlq [\phi]_{\ssl}^2}-\tl=
\frac{1}{\ml}\dfrac{\ino |\nabla \phi|^2 }{<[\phi]_{\ssl}^2>_{\ssl}}-\tl=
$$
$$
\frac{1}{\ml}\dfrac{\ino |\nabla \phi|^2 }{<\phi^2>_{\ssl}-<\phi>^2_{\ssl}}-\tl\geq
\frac{1}{\ml}\dfrac{\ino |\nabla \phi|^2 }{<\phi^2>_{\ssl}}-\tl=\dfrac{\ino |\nabla \phi|^2 }{\ino \rlq \phi^2}-\tl=
\dfrac{\ino |\nabla \phi|^2 - \tl \ino \rlq \phi^2 }{\ino \rlq \phi^2}
$$

and for $\phi \in H^{1}_0(\om)$ the equality holds if and only if $<\phi>_{\ssl}=0$. Therefore, if the equality
$\nu_{1,\lm}=\sg_{1}(\all,\pl)$ would hold, then any
eigenfunction $\phi_1$ of $\nu_{1,\lm}$ would satisfy $<\phi_1>_{\ssl}=0$, which is obviously impossible since $\phi_1$
must have constant sign. Therefore, \rife{xi1s} holds as well.

\bigskip

Concerning $D_{\al,\psi}\Phi(\lm,\al,\psi)$ we have,

\bpr\label{pr2.2} For any positive solution $(\all,\pl)$ of {\rm$\prl$} with $\lm\geq 0$, the kernel of
$D_{\al,\psi}\Phi(\lm,\all,\pl)$ is empty if and only if the equation, {\rm
\beq\label{contr}
-\Delta \phi-\tl \rlq [\phi]_{\ssl}=0,\quad \phi\in C^{2,r}_0(\ov{\om}\,)
\eeq}
admits only the trivial solution, or equivalently, if and only if $0$ is not an eigenvalue of $L_{\ssl}$.
\epr
\proof
If $\phi\in H^{1}_0(\om)$ solves \rife{contr} and since $\om$ is of class $C^3$, then by standard elliptic regularity theory
and a bootstrap argument we have
$\phi\in C^{2,r}_0(\ov{\om}\,)$. Therefore, in particular $0$ is not an eigenvalue of $L_{\ssl}$ if and only if \rife{contr}
admits only the trivial solution.\\
Suppose first that there exists a non-vanishing pair $(s,\phi)\in \R \times C^{2,r}_0(\ov{\om}\,) $ such that
$$D_{\al,\psi}\Phi(\lm,\all,\pl)[s,\phi]=(0,0).$$
Then the equation $\left.\ino \left(D_\psi \rla[\phi]+d_{\al}\rla [s]\right)\right|_{(\al,\psi)=(\all,\pl)}=0$ takes the form,
$$
p \ino \left(\lm\rlq \phi+\rlq s\right)=0,
$$
which is equivalent to $s=s_{\ssl}=- \lm<\phi>_{\ssl}$. Substituting this relation into the first equation,
$L_{\ssl}[\phi]=D_\psi F(\lm,\all,\pl)[\phi]+d_{\al}F(\lm,\all,\pl) [s_{\ssl}]=0$,
we conclude that $\phi$ is a non-trivial, classical
solution of \rife{contr}.\\ This shows one part of the claim, while on the other side, if a non-trivial, classical
solution of \rife{contr} exists, then by arguing in the other way around, obviously we can find some $(s,\phi)\neq (0,0)$
such that $D_{\al,\psi}\Phi(\lm,\all,\pl)[s,\phi]=(0,0)$, as claimed.
\finedim

\bigskip

A relevant identity is satisfied by any eigenfunction which we summarize in the following,
\ble\label{lstr}
Let $(\all,\pl)$ be a solution of {\rm $\prl$} and let $\phi_k$ be any eigenfunction of an eigenvalue $\sg_k=\sg_k(\all,\pl)$.
Then the following identity holds,
\beq\label{2907.6n}
\frac{1}{\ml}<\phi_k>_{\ssl}\equiv(\all+\lm<\pl>_{\ssl})<\phi_k>_{\ssl}=(\lm(p-1)+\sg_k)<\pl [\phi_k]_{\ssl}>_{\ssl}.
\eeq
\ele
\proof
The left hand side identity in \rife{2907.6n} is an immediate consequence of the following identity,
$$
(\all+\lm<\pl>_{\ssl})= \frac{1}{\ml}\left(\ino \rlq (\all+\lm\pl)\right)\equiv\frac{1}{\ml}.
$$
Therefore, we just need to prove the second equality.
By assumption $\phi_k$ satisfies \rife{lineqm} which we multiply by $\pl$ and integrate by parts to obtain,
$$
\ino \rl\phi_k=(\tl+\sg_k)\ino \rlq \pl[\phi_k]_{\ssl}.
$$

Dividing by $\ml$ and since $\rl=\rlq(\all+\lm\pl)$ we find that,

$$
\all<\phi_k>_{\ssl}+\lm <\pl\phi_k>_{\ssl}=(\tl+\sg_k) <\pl[\phi_k]_{\ssl}>_{\ssl}.
$$

The conclusion immediately follows by observing that
$$<\pl\phi_k>_{\ssl}=<\pl[\phi_k]>_{\ssl}+<\pl>_{\ssl}<\phi_k>_{\ssl}.$$
\finedim

\bigskip

We state now the result needed to describe branches of solutions
of $\prl$ at regular points.
\ble\label{lem1.1} Let $(\al_{\sscp \lm_0},\psi_{\sscp \lm_0})$ be a positive solution of {\rm $\prl$} with $\lm=\lm_0\geq 0$.\\
If $0$ is not an eigenvalue of $L_{\sscp \lm_0}$, then:\\
$(i)$ $D_{\al,\psi}\Phi(\lm_0, \al_{\sscp \lm_0}, \psi_{\sscp \lm_0})$ is an isomorphism;\\
$(ii)$ There exists an open neighborhood $\mathcal{U}\subset A_{\sscp \om}$ of $(\lm_0,\al_{\sscp \lm_0},\psi_{\sscp \lm_0})$ such
that the set of solutions of
{\rm $\prl$} in $\mathcal{U}$ is a real analytic curve of positive solutions $J\ni\lm\mapsto (\all,\pl)\in B$, for
suitable neighborhoods $J$ of $\lm_0$ and $B$ of $(\al_{\sscp \lm_0},\psi_{\sscp \lm_0})$ in
$(0,+\ii)\times C^{2,r}_{0,+}(\ov{\om}\,)$.\\
$(iii)$ In particular if $(\al_{\sscp \lm_0},\psi_{\sscp \lm_0})=(\al_{\sscp 0},\psi_{\sscp 0})=(1,G[1])$, then
$(\all,\pl)=(1,\psi_{\sscp 0})+\mbox{\rm O}(\lm)$ as $\lm\to 0$.
\ele

\proof
For the sake of simplicity, in the rest of this proof we set $\lm_0=\lm$.\\
By the construction of $A_{\sscp \om}$, the map
$F$ as defined in \rife{eF} is jointly analytic in a suitable neighborhood of $(\lm,\all,\pl)$. As a consequence, whenever $(i)$ holds,
then $(ii)$ is an immediate consequence of the analytic implicit function theorem,
see for example Theorem 4.5.4 in \cite{but}. In particular $(iii)$ is a straightforward consequence of $(ii)$. Therefore, we are
just left with the proof of $(i)$.\\

$(i)$  Let $(t,f)\in \R\times C^{r}(\ov{\om}\,)$, then we have to prove that the equation,
$$
D_{\al,\psi}\Phi(\lm,\all,\pl)[s,\phi] =(t,f),
$$
admits a unique solution $(s,\phi)\in \R\times C^{2,r}_0(\ov{\om}\,)$. In view of Proposition \ref{pr2.2},
it is enough to prove that it actually admits just
one solution. The equation
$$
\left. \ino \left(D_\psi \rla[\phi]+d_{\al}\rla [s]\right)\right|_{(\al,\psi)=(\all,\pl)}=t
$$
takes the form,
$$
p \ino \left(\lm\rlq \phi+\rlq s\right)=t,
$$
which we solve as follows,
$$
s=s_{\ssl}=\frac{t}{p\ml}- \lm<\phi>_{\ssl}.
$$
Substituting this relation into the first equation, that is,
$$
D_\psi F(\lm,\all,\pl)[\phi]+d_{\al}F(\lm,\all,\pl) [s_{\ssl}]=f,
$$
by standard elliptic estimates we just need to find a solution $\phi_{\sscp f, \sscp t}\in C^{0}(\ov{\om}\,)$ of the equation,
\beq\label{un1}
\phi=G[\tl\rlq[\phi]_{\ssl}]+G\left[ f_t\right], \quad \mbox{ where }f_t=f+t\frac{\rlq}{\ml}.
\eeq
Let $Y_{0,r}=Y_0\cap C^{r}(\ov{\om}\,)$ (see \rife{y0}) and $Y^{\perp}_{0,r}$  be the orthogonal complement of
$Y_{0,r}$ in $C^r(\ov{\om}\,)$.
Since $Y^{\perp}_{0,r} \oplus Y_{0,r} = C^{r}(\ov{\om}\,)$, then after projection on $Y^{\perp}_{0,r}$ and $Y_{0,r}$, we see that
\rife{un1} is equivalent to the system,
$$
\graf{[\phi]_{\ssl}= [\;G [\tl\rlq[\phi]_{\ssl}]+G[f_t] \;]_{\ssl},\\ \\<\phi>_{\ssl}=<G [\tl\rlq[\phi]_{\ssl}]+G[f_t]>_{\ssl}.}
$$
Since $0\notin \{\sg_j(\all,\pl)\}_{j\in\N}$, then $I-T_{0}$ (see \rife{T0}) is an isomorphism of $Y_{0,r}$ onto itself.
Therefore, the first equation, which has the form $(I-T_0)([\phi]_{\ssl})=[G[f_t]]_{\ssl}$,
has a unique solution $\varphi_{\sscp f, \sscp t}\in Y_{0,r}$. Let
$$
\phi_{\sscp f, \sscp t}:=G[\tl\rlq \varphi_{\sscp f, \sscp t}]+G[f_t]=\varphi_{\sscp f, \sscp t} +<G[\tl\rlq \varphi_{\sscp f, \sscp t}]+G[f_t]>_{\ssl}.
$$
By standard elliptic regularity theory we find that $\phi_{\sscp f, \sscp t}\in C^{2,r}_0(\ov{\om}\,)$, which in particular
is the unique solution of both equations.\finedim

\bigskip

The next result is about the transversality condition needed in section \ref{sec6}.
\bte\label{spectral1} Let $(\all,\pl)$ be a positive solution of {\rm$\prl$}. Suppose that any eigenfunction $\phi_k$ of a fixed
vanishing eigenvalue $\sg_k=\sg_k(\all,\pl)=0$ satisfies $<\phi_k>_{\ssl}\neq 0$. Then $\sg_k$ is simple, that is,
it admits only one eigenfunction.\\
Let $\om\subset \R^2$ be symmetric and convex with respect to the coordinate directions $x_i$, $i=1,2$ and let $(\all,\pl)$ be a positive solution of {\rm$\prl$} with $\lm>0$.
Suppose that $\sg_k=\sg_k(\all,\pl)=0$ and let $\phi_k$ be any
corresponding eigenfunction. Then:\\
$(i)$ $<\phi_k>_{\ssl}\neq 0$;\\
$(ii)$ $\sg_k(\all,\pl)$ is simple, that is, it admits at most one eigenfunction.
\ete
\proof
If any eigenfunction $\phi_k$ of a vanishing eigenvalue $\sg_k=0$ satisfies $<\phi_{k}>_{\ssl}\neq 0$, then
there can be at most one such eigenfunction. Indeed, if there were more than one such eigenfunctions,
say $\phi_{k,\ell}$, $\ell=1,2$, then, putting
$a_\ell:=<\phi_{k,\ell}>_{\ssl}\neq 0$, $\ell=1,2$, we would find that
$\phi=\phi_{k,1}-\frac{a_1}{a_2}\phi_{k,2}$ would be an eigenfunction of
$\sg_k$ satisfying $<\phi>_{\ssl}=0$, which is a contradiction.\\
Concerning $(i)$ we argue by contradiction and assume that $<\phi_k>_{\ssl}= 0$.
In view of \eqref{lineqm} and standard elliptic regularity theory $\phi_k$ is a classical solution of,
\beq\label{0107.1}
\graf{-\Delta \phi_k =\tl \rlq \phi_k \quad \mbox{in}\;\;\om,\\ \\
\phi_k=0 \quad \mbox{on}\;\;\pa\om,
}
\eeq
satisfying $\ml <\phi_k>_{\ssl}=\ino \rlq \phi_k=0$.
However, by Theorem 3.1 in \cite{DGP} (which relies on the symmetry and convexity properties of $\om$) the
nodal line of any solution of \rife{0107.1} cannot intersect the boundary. Here one is also using the fact that $\pl\geq 0$ is
a solution of $\prl$ and that
$\rl=f_{\ssl}(\pl)$, where $f_{\ssl}:[0,+\ii)\to [0,+\ii)$ is a $C^1$ function and  $f_{\ssl}(0)\geq 0$.
Therefore, $\phi_k$ has a fixed sign in a small enough neighborhood of the boundary and then, in particular,
by the strong maximum principle we can assume without loss of generality that
$\pa_{\nu}\phi_k=\frac{\pa \phi_k}{\pa \nu}<0$ on $\pa \om$,
where $\nu$ denotes the exterior unit normal.
On the other side, by \rife{0107.1} and $<\phi_k>_{\ssl}= 0$, we see that $\inpo \pa_{\nu}\phi_k=0$,
which is the desired contradiction.\\
At this point $(ii)$ follows immediately from $(i)$ and the first part of the claim.
\finedim

\bigskip
\bigskip

The following result ensures that the first eigenvalue $\sg_1(\all,\pl)$ of a positive
variational solution must be non-negative.

\ble\label{prmin}
Let $(\all,\pl)$ be a positive variational solution of {\rm $\prl$}. Then $\sg_1(\all,\pl)\geq 0$.
\ele
\proof Since $\rl$ is a minimizer of the variational problem {\bf (VP)} in Appendix \ref{appD}, then the Taylor formula shows that,
$$
0\leq J_{\ssl}(\rl+\eps f)-J_{\ssl}(\rl)=\frac{\eps^2}{2}\mathcal{D}^{2}J_\lm(\rl)[f,f\,]+\mbox{\rm o}(\eps^2),
$$
for any $f \in C^{2}_0(\ov{\om}\,)$ such that $\ino f=0$. Therefore, for any such $f$, we find,
$$
\mathcal{D}^{2}J_\lm(\rl)[f,f\,]=\frac1p\ino (\rl)^{\frac1p-1}f^2-\lm\ino fG[f]\geq 0,
$$
which is well defined since $\rl\geq \all^p>0$ on $\ov{\om}$. Clearly, for any $\varphi\in C^{2}_0(\ov{\om}\,)$
such that $\ino \rlq \varphi=0$ we can choose $f=\rlq \varphi$ and deduce that,
$$
\ino (\rl)^{\frac1p-1}f^2-p\lm\ino fG[f]=\ino \rlq\varphi^2-p\lm\ino \rlq \varphi \,G[\rlq \varphi]\geq 0.
$$
Therefore, we have
$$
\mathcal{A}(\phi):=\ino \rlq [\phi]_{\ssl}^2-\tl\ino \rlq  [\phi]_{\ssl} G[\rlq [\phi]_{\ssl}]\geq 0,\,\forall\,\phi\in C^{2}_0(\ov{\om}\,),
$$
where we recall that $[\phi]_{\ssl}=\phi\,-<\phi>_{\ssl}$ and $\tl=p\lm$.
On the other side, letting $\phi_1$ be the first eigenfunction
of \rife{lineq0} whose eigenvalue is $\sg_1$, from \rife{lineq9m} we see that,
$$
0\leq \dfrac{\mathcal{A}(\phi_1)}{\ml}=<[\phi_1]_{\ssl}^2>_{\ssl}-\tl<[\phi_1]_{\ssl} G[\rlq [\phi_1]_{\ssl}]>_{\ssl}=
$$
$$
<[\phi_1]_{\ssl}^2>_{\ssl}-\frac{\tl}{\tl+\sg_1}<[\phi_1]_{\ssl}\phi_1>_{\ssl}=<[\phi_1]_{\ssl}^2>_{\ssl}\frac{\sg_1}{\tl+\sg_1},
$$
and then from \rife{lineq2} we immediately conclude that $\sg_1\geq0$.\finedim

\bigskip
\bigskip

Finally we obtain a generalization of the bending result of \cite{CrRab} for solutions of $\prl$, which
is an improvement of Lemma \ref{lem1.1} in case of simple and vanishing eigenvalues satisfying a suitable transversality
condition. We need here to exploit the modified spectral setting in its full strength.
\bpr\label{pr3.1}
Let $(\all,\pl)$ be a positive solution of {\rm$\prl$} with $\lm>0$ and suppose that the  $k$-th
eigenvalue $\sg_k(\all,\pl)=0$ is simple, that is, it admits only one eigenfunction, $\phi_k\in C^{2,r}_0(\ov{\om}\,)$.
If $<\phi_k>_{\ssl}\neq 0$,
then there exists $\eps>0$, an open neighborhood $\mathcal{U}$ of $(\lm,\all,\pl)$ in $A_{\sscp \om}$ and
a real analytic curve $(-\eps,\eps)\ni s \mapsto (\lm(s), \al(s),\psi(s))$ such that
$(\lm(0), \al(0),\psi(0))=(\lm,\all,\pl)$ and the set of solutions of {\rm $\prl$} in $\mathcal{U}$ has the form $(\lm(s),\al(s),\psi(s))$,
where $(\al(s),\psi(s))$ is a solution of {\rm$\prl$} for $\lm=\lm(s)$ for any $s\in (-\eps,\eps)$, with
$\psi(s)=\pl+s\phi_k+\vxi(s)$, and
\beq\label{2907.0}
<[\phi_k]_{\sscp \lm(s)},\vxi(s)>_{\sscp \lm(s)}=0,\quad s\in (-\eps,\eps).
\eeq
Moreover it holds,
\beq\label{2907.1}
\vxi(0)\equiv 0\equiv \vxi^{'}(0),\quad \al^{'}(0)=-\lm <\phi_k>_{\ssl},\quad\lm^{'}(0)=0,\quad \psi^{'}(0)= \phi_k,
\eeq
and either $\lm(s)=\lm$ is constant in $(-\eps,\eps)$ or
$\lm^{'}(s)\neq 0$, $\sg_k(s)\neq 0$ in $(-\eps,\eps)\setminus\{0\}$, $\sg_k(s)$ is simple in $(-\eps,\eps)$ and

\beq\label{2907.5}
<[\phi_k]_{\ssl},\pl>_{\ssl}\neq 0\mbox{ and $<[\phi_k]_{\ssl},\pl>_{\ssl}$  has the same sign as }<\phi_k>_{\ssl},
\eeq

\beq\label{2907.10}
\dfrac{\sg_k(s)}{\lm^{'}(s)}=\dfrac{p<[\phi_k]_{\ssl},\psi_{\ssl}>_{\ssl}+\mbox{\rm o}(1)}
{<[\phi_k]_{\ssl}^2>_{\ssl}+\mbox{\rm o}(1)},\mbox{ as }s\to 0.
\eeq
\epr
\proof Since $<\phi_k>_{\ssl}\neq 0$ by assumption, then \rife{2907.5} is an immediate consequence of \rife{2907.6n}.\\
Take $\eps>0$ small enough and
$\dt_i=\dt_i(\eps)>0$, $i=1,2,3$, such that
$$
(\lm+\mu,\all+\beta,\pl+s\phi_k+\vxi)\subset A_{\sscp \om},
\;\forall\,(s,\mu,\beta,\vxi)\in (-\eps,\eps)\times (-\dt_1(\eps),\dt_1(\eps))\times (-\dt_2(\eps),\dt_2(\eps))\times B_{k,\eps}^{\perp},
$$
where
$$
B_{k,\eps}^{\perp}=\left\{\vxi\in C^{2,r}_0(\ov{\om}\,)\,:\,< [\phi_k]_{\ssl},\vxi>_{\ssl}=0,\;\|\vxi\|_{C^{2,r}_0(\ov{\om}\,)}< \dt_3(\eps)\right\}.
$$

Next, let us introduce the map,
$$
\Phi_0: (-\eps,\eps)\times (-\dt_1(\eps),\dt_1(\eps))\times (-\dt_2(\eps),\dt_2(\eps))\times B_{k,\eps}^{\perp}\to \R \times  C^{r}(\ov{\om}\,),
$$
defined as follows,
$$
\Phi_0(s,\mu,\beta,\vxi)=\left(\begin{array}{cl}-\Delta (\pl+s\phi_k+\vxi)-(\all+\beta +(\lm+\mu)(\pl+s\phi_k+\vxi))^p\\
\\-1+\ino (\all+\beta +(\lm+\mu)(\pl+s\phi_k+\vxi))^p\end{array}\right).
$$

Clearly, with the notations introduced in \rife{ePhi}, we have
$$
\Phi_0(0,0,0,0)=\Phi(\lm,\all,\pl)=(0,0).
$$
Let us define,
$$
X_{k}^{\perp}=\left\{\phi\in C^{2,r}_0(\ov{\om}\,)\,:\,< [\phi_k]_{\ssl},\phi>_{\ssl}=0\right\},
$$
then the differential of $\Phi_0$ with respect to
$(\mu,\beta,\vxi)$ at $(0,0,0,0)$ acts on a triple $(s_\mu, s_\be,\phi)\in \R\times\R\times X_{k}^{\perp}$ as follows,
$$
D_{\mu,\beta,\vxi}\Phi_0(0,0,0,0)[s_\mu,s_\be,\phi]=
\left(\begin{array}{cl} D_\psi F(\lm,\all,\pl)[\phi]-p\rlq \pl s_\mu-p\rlq s_\be \\ \\
\tl \ino \rlq \phi +p s_{\mu}\ino \rlq \pl +p s_{\be}\ino \rlq \end{array}\right),
$$
with $\tl=p\lm$ and $D_\psi F$ as in \rife{lin}.
The crux of the argument is to prove the following:
\ble\label{lem.crux}
$D_{\mu,\beta,\vxi}\Phi_0(0,0,0,0)$ is an isomorphism of $\R\times\R\times X_{k}^{\perp}$ onto $\R \times  C^{r}(\ov{\om}\,)$.
\ele
\proof
We will prove that for each $(t,f)\in \R \times  C^{r}(\ov{\om})$, the vectorial equation
$$
D_{\mu,\beta,\vxi}\Phi_0(0,0,0,0)[s_\mu,s_\be,\phi]=\left(\begin{array}{cl}f\\
\\t\end{array}\right),
$$
admits a unique solution $(s_\mu,s_\be,\phi)\in \R\times\R\times X_{k}^{\perp}$. From the second equation we deduce that,
$$
s_{\be}=-\frac{t }{p\ml}-<\pl>_{\ssl} s_\mu -\lm  <\phi>_{\ssl},
$$
where we recall that $\ml=\ino\rlq$.
Therefore, by substituting into the first equation, we find that the pair $(s_\mu,\phi)$ solves the equation,
\beq\label{26.7.1}
-\Delta \phi-\tl \rlq [\phi]_{\ssl}-p\rlq [\pl]_{\ssl} s_\mu = f+ t \frac{\rlq}{\ml}.
\eeq

Let us write $C^r(\ov{\om}\,)=Y_k \oplus R$, where $Y_k =\mbox{span}\{\phi_k\}$ and
$R=\left\{f\in C^{r}(\ov{\om}\,)\,:\,\ino \phi_k f =0 \right\}$.
Then, projecting \rife{26.7.1} along $Y_k$, we find that,

\beq\label{26.7.2}
-m_{\ssl}p <[\pl]_{\ssl}, \phi_k>s_\mu = \ino\left(f+ t \frac{\rlq}{\ml}\right)\phi_k.
\eeq

At this point it follows from \rife{2907.5} that \rife{26.7.2} admits a unique solution $s_\mu$.
So we are left with showing that the projection of \rife{26.7.1} onto $R$ admits
a unique solution as well.

First of all, let us observe that, since $\phi_k$ satisfies \eqref{lineqm} with $\sg_k=0$,
then the projection of \rife{26.7.1} onto $R$ takes the form,
$$
-\Delta \phi-\tl \rlq [\phi]_{\ssl}= g,
$$
where
$$
g=P_R\left(f+  t \frac{\rlq}{\ml}+p\rlq [\pl]_{\ssl} s_\mu\right),
$$
where $P_R$ denotes the projection operator. Since $\all> 0$ and $\lm>0$ by assumption, then $\rl\geq \all^p$ in $\ov{\om}$
and then, by standard elliptic estimates, it is enough to show that the equation,
\beq\label{26.7.4}
 \phi=\tl G[\rlq [ \phi]_{\ssl}] +G[\rlq g_1],
\eeq
admits a unique solution $ \phi\in Y_{k}^{\perp}:=\{f \in C^{r}(\ov{\om}\,)\,:\,<[\phi_k]_{\ssl},f>_{\ssl}=0\}$,
for any fixed $g_1\in C^r(\ov{\om})$ satisfying $< g_1,\phi_k>_{\ssl}=0$,
where $<[\phi_k]_{\ssl}^2>_{\ssl}=1$. Next, let us decompose,
$$
\phi=c_0+c_k[\phi_k]_{\ssl}+\varphi,\quad g_1=b_0+b_k[\phi_k]_{\ssl}+g_{2},
$$
where $c_0=<\phi>_{\ssl}$, $b_0=<g_1>_{\ssl}$, $c_k=<[\phi_k]_{\ssl},\phi>_{\ssl}$, $b_k=< [\phi_k]_{\ssl},g_1>_{\ssl}$,
$\{\varphi,g_2\}\in Y_{k,0}:=Y_{0}\cap Y_{k}^{\perp}$, $Y_{0}=\{f \in C^{r}(\ov{\om})\,:\,<f>_{\ssl}=0\}$.
It is easy to check, using \rife{lineq9m}, that we have $c_k=0$, $b_k+b_0<\phi_k>_{\ssl}=0$ and
$$
<G[\rlq \varphi],[\phi_k]_{\ssl}>_{\ssl}=0=<G[\rlq g_2],[\phi_k]_{\ssl}>_{\ssl}.
$$
At this point, by using also \rife{lineq9m}, then a lengthy
evaluation shows that \rife{26.7.4} is equivalent to the system,
$$
\graf{c_0= b_0 <G[\rlq]>_{\ssl}+\tl <G[\rlq \varphi]>_{\ssl}+\frac{b_k}{\tl}<\phi_k>_{\ssl}+<G[\rlq g_2]>_{\ssl},\\ \\
\varphi-\tl [G[\rlq \varphi]]_{\ssl}= [G[\rlq g_2]]_{\ssl} +b_0 \left([G[\rlq]]_{\ssl}-
<[G[\rlq]]_{\ssl},[\phi_k]_{\ssl}>_{\ssl}[\phi_k]_{\ssl}\right),\\ \\
b_0<[G[\rlq]]_{\ssl},[\phi_k]_{\ssl}>_{\ssl} +\frac{b_k}{\tl}=0,}
$$
with the constraint $b_k+b_0<\phi_k>_{\ssl}=0$.
By using  \rife{lineq9m} once more we see that $<[G[\rlq]]_{\ssl},[\phi_k]_{\ssl}>_{\ssl}=\frac{<\phi_k>_{\ssl}}{\tl}$,
whence the third equation is
equivalent to $b_k+b_0<\phi_k>_{\ssl}=0$ and we are left with the first two equations. However the second equation takes the form,
$$
\left.(I-T_0)\right|_{Y_{k,0}}(\varphi)=f_2,
$$
for a suitable $f_2\in Y_{k,0}$, where $T_0$ is given in \rife{T0}. Since $\left.(I-T_0)\right|_{Y_{k,0}}(h)=0$ has only the trivial solution,
then by the Fredholm alternative it
admits a unique solution, say $\varphi_2$. As a consequence $\varphi_2$ uniquely defines $c_0$ in the first equation, as claimed.
\finedim

\bigskip
\bigskip

In view of Lemma \ref{lem.crux}, the existence of the real analytic map $(-\eps,\eps)\ni s \mapsto (\lm(s), \al(s),\psi(s))$
satisfying the desired properties (including \rife{2907.0}) follows by the analytic implicit function theorem (\cite{but}).
Therefore, we are left with the proof of \rife{2907.1} and \rife{2907.10}.\\
First of all, differentiating the constraint in $\prl$ with $\lm(s)=\lm+\mu(s)$, $\al(s)=\all+\be(s)$, $\psi(s)=\pl+s\phi_k+\vxi(s)$ we obtain,
$$
p\ino\rlq\left(\beta^{'}(0)+\mu^{'}(0)\pl+\lm (\phi_k+\vxi^{'}(0))\right)=0,
$$
which immediately implies that,
\beq\label{2907.2}
\beta^{'}(0)=-\mu^{'}(0)<\pl>_{\ssl}-\lm<\phi_k>_{\ssl}-\lm<\vxi^{'}(0)>_{\ssl}.
\eeq
Thus, by differentiating the equation in $\prl$ we obtain,
$$
-\Delta (\phi_k+\vxi^{'}(0))=p\rlq\left(\beta^{'}(0)+\mu^{'}(0)\pl+\lm (\phi_k+\vxi^{'}(0)) \right),
$$
and so, by using \rife{lineqm} and \rife{2907.2}, we conclude that,
\beq\label{2907.3}
-\Delta \vxi^{'}(0)=\tl \rlq [\vxi^{'}(0)]_{\ssl}+p\mu^{'}(0)\rlq [\pl]_{\ssl},
\eeq
where, since $\vxi(0)\equiv 0$, then differentiating \rife{2907.0}, we also have $<[\phi_k]_{\ssl},\vxi^{'}(0)>_{\ssl}=0$. Therefore, multiplying \rife{2907.3} by
$\phi_k$ and integrating by parts, we find that,
$$
0=\ml\tl<[\phi_k]_{\ssl},\vxi^{'}(0)>_{\ssl}=\ino (-\Delta \phi_k)\vxi^{'}(0)=
$$
$$
\ino \phi_k(-\Delta \vxi^{'}(0))=\ino\tl \rlq [\vxi^{'}(0)]_{\ssl}\phi_k+
p\mu^{'}(0)\ino \rlq [\pl]_{\ssl}\phi_k=
$$
$$
\ml\tl<[\phi_k]_{\ssl},\vxi^{'}(0)>_{\ssl}+\ml p \mu^{'}(0)<[\phi_k]_{\ssl},\pl>_{\ssl},
$$
that is,
$$
\mu^{'}(0)<[\phi_k]_{\ssl},\pl>_{\ssl}=0,
$$
which, in view of \rife{2907.5}, yields $\mu^{'}(0)=0$. Inserting this information back in \rife{2907.3},
in view of $<[\phi_k]_{\ssl},\vxi^{'}(0)>_{\ssl}=0$, and since
we assumed that $\sg_k$ is simple, then we conclude that $\vxi^{'}(0)\equiv 0$. Therefore, \rife{2907.2} shows that
$\beta^{'}(0)=-\lm<\phi_k>_{\ssl}$, which concludes the proof of \rife{2907.1}.\\
Concerning \rife{2907.10}, let us first observe that, since $\lm$ is real analytic in $(-\eps,\eps)$, then either
$\lm(s)=\lm$ is constant or $\lm^{'}(s)\neq 0$ in $(-\eps,\eps)\setminus\{0\}$ for $\eps$ sufficiently small. Also,
$\sg_k(s)=\sg_k(\al(s),\psi(s))$ is simple in $(-\eps,\eps)$
by the analytic perturbation theory of simple eigenvalues, see Proposition 4.5.8 in \cite{but}.
Clearly, in view of
\rife{2907.5}, we can also assume that $<[\psi(s)]_{\sscp \lm(s)},\phi_k(s)>_{\sscp \lm(s)}\neq 0$ for any such $s$.
Next, by arguing as in \rife{2907.3}
we see that,
\beq\label{2907.12}
-\Delta \psi^{'}(s)=p\lm(s) (\rh_{\sscp \lm(s)})^{\frac1q} [\psi^{'}(s)]_{\sscp \lm(s)}+p\lm^{'}(s)(\rh_{\sscp \lm(s)})^{\frac1q} [\psi(s)]_{\sscp \lm(s)},
\eeq
and so, multiplying \rife{2907.12} by $\phi_k(s)$,  integrating by parts and using \rife{lineqm} with $\lm=\lm(s)$, then we conclude that,
\beq\label{2907.13}
\sg_k(s)\ino (\rh_{\sscp \lm(s)})^{\frac1q} [\phi_k(s)]_{\sscp \lm(s)}\psi^{'}(s)=p\lm^{'}(s)\ino (\rh_{\sscp \lm(s)})^{\frac1q} [\phi_k(s)]_{\sscp \lm(s)}\psi(s).
\eeq

At this point we observe that, since  $<[\psi(s)]_{\sscp \lm(s)},\phi_k(s)>_{\sscp \lm(s)}\neq 0$, and
$\ino \rh_{\sscp \lm(s)} [\phi_k(s)]_{\sscp \lm(s)}\psi^{'}(s)=\mbox{\rm o}(1)+\ino [\phi_k]_{\ssl}^2$, as $s\to 0$, then,
if $\lm^{'}(s)\neq 0$, we can write,
$$
\dfrac{\sg_k(s)}{\lm^{'}(s)}=\dfrac{p\ino (\rh_{\sscp \lm(s)})^{\frac1q} [\phi_k(s)]_{\sscp \lm(s)}\psi(s)+\mbox{\rm o}(1)}{\ino (\rh_{\sscp \lm(s)})^{\frac1q} [\phi_k(s)]_{\sscp \lm(s)}\psi^{'}(s)+\mbox{\rm o}(1)},
$$
for any $s$ small enough, thus,
$$
\dfrac{\sg_k(s)}{\lm^{'}(s)}=\dfrac{p<[\phi_k]_{\ssl},\psi_{\ssl}>_{\ssl}+\mbox{\rm o}(1)}
{<[\phi_k]_{\ssl}^2>_{\ssl}+\mbox{\rm o}(1)},\mbox{ as }s\to 0,
$$
which is \rife{2907.10}.
\finedim

\bigskip

\section{\bf Monotonicity of solutions}\label{sec3}
In this section we are concerned with the monotonicity of the energy
$$
2 E_{\ssl}=\ino \rl \pl,
$$
and of $\all$ for positive solutions of $\prl$ with $\sg_1(\all,\pl)>0$, see \rife{4.1}. The latter property is detected by making use of the spectral setting set up in the previous section. We recall that $\ml=\ino \rlq$ and $L_{\ssl}$
 is the linearized operator defined in \rife{eLl}.
\bpr\label{pr-enrg} Let $(\al_{\ssl_0},\psi_{\ssl_0})$ be a positive solution of {\rm $\prl$} with $\lm_0\geq 0$
and suppose that $0$ is not an eigenvalue of $L_{\ssl_0}$.
Then, locally near $\lm_0$, the map $\lm\mapsto (\all,\pl)$ is a real analytic simple curve of positive solutions and
{\rm
$$
2 E_{\ssl}=\ino \rl \pl\quad \mbox{\it and }\quad \vl=\dfrac{d \pl}{d \lm},
$$}

are real analytic functions of $\lm$ with $\vl \in C^{2}_0(\ov{\om}\,)$. In particular we have:
{\rm
\beq\label{8.12.11}
\frac{d E_{\ssl}}{d \lm} = \ino \rl \vl=\ml \tl<[\vl]_{\ssl}, [\pl]_{\ssl}>_{\ssl}+\ml p\|[\pl]_{\ssl}\|^2_{\ssl}.
\eeq}

Finally, if $\sg_1=\sg_1(\all,\pl)> 0$, then,
{\rm
\beq\label{21.06.17}
<[\vl]_{\ssl}, [\pl]_{\ssl}>_{\ssl}\geq\frac{\sg_1}{p} \|[\vl]_{\ssl}\|^2_{\ssl},
\eeq
}
and in particular $\frac{d \el}{d\lm}>0$.
\epr

\proof Clearly, $\frac{d \el}{d\lm}>0$ will follow immediately once we have \rife{8.12.11}, \rife{21.06.17}.\\
Since by assumption $0$ is not an eigenvalue of $L_{\ssl}$, then we can apply Lemma \ref{lem1.1}. Therefore, locally around $(\al_{\ssl_0},\psi_{\ssl_0})$,
$(\all,\pl)$ is a real analytic function of $\lm$ and then, by standard elliptic estimates, we see that
$\vl\in C^{2}_0(\ov{\om}\,)$ is a classical solution of,
$$
-\Delta \vl =\tl \rlq \vl +p \rlq\pl +p\rlq \frac{d\all}{d\lm},
$$
where $\frac{d\all}{d\lm}$ can be computed by the unit mass constraint in $\prl$, that is
$$
p\frac{d\all}{d\lm}=-\tl  <\vl>_{\ssl}-p<\pl>_{\ssl}.
$$
Therefore, we conclude that $\vl$ is a solution of,
\beq\label{1b1}
-\Delta \vl =\tl \rlq [\vl]_{\ssl} +p \rlq[\pl]_{\ssl}.
\eeq
By using $\prl$, we also have $E_{\ssl}=\frac12\ino |\nabla \pl|^2$, and in particular it holds,
$$
\frac{d}{d \lm}E_{\lm} =\ino (\nabla \vl,\nabla \pl)=-\ino \vl (\Delta \pl)=\ino \rl\vl,
$$
which proves the first equality sign in \rife{8.12.11}. At this point, by using $\prl$ and \eqref{1b1}, we see that,
\beq\label{2b1}
\ino\rl \vl=\ino -(\Delta \pl)\vl=\ino -\pl(\Delta \vl)= \ml \tl <\pl, [\vl]_{\ssl}>_{\ssl} +\ml  p <\pl, [\pl]_{\ssl}>_{\ssl},
\eeq
which proves the second equality sign in \rife{8.12.11}.\\

For the last part of the statement, let
$$
[\pl]_{\ssl}=\sum\limits_{j=1}^{+\ii}\xi_j[\phi_j]_{\ssl}, \quad
[\vl]_{\ssl}=\sum\limits_{j=1}^{+\ii}\beta_j[\phi_j]_{\ssl},
$$
$$
\xi_j=\ino \rlq[\phi_j]_{\ssl}[\pl]_{\ssl}, \quad
\beta_j=\ino \rlq [\phi_j]_{\ssl}[\vl]_{\ssl},
$$
be the Fourier expansions  of $[\pl]_{\ssl}$ and $[\vl]_{\ssl}$ in $Y_0$ (see \rife{y0}), with
respect to the normalized projections $[\phi_j]_{\ssl}$, satisfying $\|[\phi_j]_{\ssl}\|_{\ssl}=1$.
After multiplying \rife{1b1} by $\phi_j$, using \rife{lineq0} and integrating by parts, we have,
\beq\label{lamq31}
\sg_{j}\ino\rlq[\phi_j]_{\ssl}[\vl]_{\ssl} =p
\ino\rlq [\phi_j]_{\ssl}[\pl]_{\ssl},\mbox{ that is }\sg_j\beta_j=p{\xi_j},
\eeq
where $\sg_{j}=\sg_{j}(\all,\pl)$. As a consequence, since $\sum\limits_{j=1}^{+\ii} (\beta_j)^2=<[\vl]_{\ssl}^2>_{\lm}$ and
$\sum\limits_{j=1}^{+\ii} (\xi_j)^2=<[\pl]_{\ssl}^2>_{\lm}$, then we find that,
\beq\label{event1}
p<[\pl]_{\ssl},[\vl]_{\ssl}>_{\ssl}=p\sum\limits_{j=1}^{+\ii} \xi_j\beta_j=
\sum\limits_{j=1}^{+\ii} \sg_{j}(\beta_j)^2\geq \sg_{1}<[\vl]_{\ssl}^2>_{\lm},
\eeq
which proves \rife{21.06.17}.
\finedim

\bigskip

Next we prove the monotonicity of $\all$ whenever $\sg_1(\all,\pl)>0$.

For this and later purposes it is convenient to use the auxiliary function $\xil=\lm \pl$ which satisfies
\beq\label{ul.1}
\graf{-\Delta \xil =\lm \left(\all+\xil\right)^p\quad \mbox{in}\;\;\om\\ \\
\bigintss\limits_{\om}  \left(\all+\xil\right)^p=1\\ \\
\xil\geq0 \;\; \mbox{in}\;\;\om, \quad \xil=0 \;\; \mbox{on}\;\;\pa\om \\ \\
\all\geq0
}
\eeq
where $\rl=\left(\all+\xil\right)^p$, $\rlq=\left(\all+\xil\right)^{p-1}$.

\bpr\label{pr3.2.best}
Let $(\all,\pl)$ be a positive solution of {\rm $\prl$} {\rm (}or either $(\all,\xil)$ a solution of {\rm \rife{ul.1}}, $\all>0${\rm )}
 with $\lm\geq 0$.
If $\sg_1=\sg_1(\all,\pl)>0$ then $\dfrac{d \all}{d \lm}<0$.
\epr
\proof
We recall that, by Lemma \ref{lem1.1}, if $0$ is not an eigenvalue of $L_{\ssl}$, then $\pl$ and $\all$
are locally real analytic functions of $\lm$. Therefore, $\xil$ is also real analytic as a function of $\lm$ and
$\wl=\frac{d\xil}{d\lm}\in C^{2}_0(\ov{\om}\,)$ satisfies
\beq\label{wl}
\graf{-\Delta \wl =\tl \rlq [\wl]_{\ssl}+\rl\quad \mbox{in}\;\;\om\\ \\
\wl=0 \quad \mbox{on}\;\;\pa\om
}
\eeq
where, since $\ino (\all + {\xil})^p=1$, then,
\beq\label{malfa.1}
\frac{d\all}{d\lm}=-<\wl>_{\ssl}.
\eeq
We first prove the claim for $p>1$ and $\lm>0$.
Multiplying the equation in \rife{wl} by $\wl$ and integrating by parts
we find that,
\beq\label{alm1}
\all\ino \rlq \wl+\ino \rlq \wl\xil=\ino \rl \wl=\ino |\nabla \wl|^2-\tl \ino \rlq[\wl]^2_{\ssl} \geq
\eeq
$$
\sg_{1}(\all,\pl)\ino \rlq[\wl]^2_{\ssl}.
$$

On the other side, multiplying the equation in \rife{ul.1} by $\wl$ and integrating by parts
we also find that,

\beq\label{030520.1}
\tl \ino \rlq [\wl]\xil+\ino \rl \xil=\lm \ino \rl \wl=\lm \all \ino \rlq \wl+\lm \ino \rlq \wl{\xil},
\eeq
which, after a straightforward evaluation yields,

$$
\ino\rlq \wl\xil=\frac{1}{p-1}(p<\xil>_{\ssl}+\all)\ino\rlq\wl -\frac{1}{\lm(p-1)}\ino \rl \xil.
$$
Substituting this expression of $\ino\rlq \wl\xil$ in \rife{alm1} we obtain,

$$
\tl \left(\all + {<\xil>_{\ssl}}\right)\ino \rlq \wl-\ino \rl \xil
\geq \lm(p-1)\sg_{1}(\all,\pl)\ino \rlq[\wl]^2_{\ssl}.
$$

In other words we conclude that,

$$
\tl<\wl>_{\ssl}=\frac{\tl}{\ml} \ino \rlq \wl\geq {\lm}{(p-1)}\sg_{1}(\all,\pl)\ino \rlq[\wl]^2_{\ssl}+\ino \rl \xil,
$$
where we are using the identity

$$
\all + {<\xil>_{\ssl}}=\frac{\ino \rlq(\all + {\xil})}{\ml}=\frac{\ino \rl}{\ml}=\frac{1}{\ml}.
$$

At this point, since $\xil\geq 0$,  $\frac{d\all}{d\lm}<0$ immediately follows from \rife{malfa.1}.\\
For $p=1$ and $\lm>0$ after a straightforward evaluation we deduce from \rife{030520.1} that $<\wl>_{\ssl}=2\el$ and then the conclusion
follows from \rife{malfa.1}. For $p\geq1$ and $\lm=0$ we have $\wl=\psi_{\sscp 0}$ and $<\wl>_{\ssl}=2E_{0}$ and the conclusion follows again from \rife{malfa.1}.
\finedim

\bigskip
\bigskip

\section{\bf Uniqueness of solutions}\label{sec4}

The following Lemma is the starting point of the proof of Theorem \ref{thmLE} and illustrates an interesting
property of $\prl$.
\ble\label{thm05} Let $p\in [1,p_{N})$.\\
$(i)$ If there exists $\ov{\lm}>0$ such that for any solution $(\all,\pl)$ with $\lm\leq\ov{\lm}$ it holds
$\sg_1(\all,\pl)>0$, then {\rm $\prl$} has at most one positive solution for any $\lm\leq\ov{\lm}$.

\medskip

$(ii)$ If there exists $\ov{\lm}>0$ such that for any solution $(\all,\pl)$ with $\lm\leq\ov{\lm}$ it holds
$\sg_1(\all,\pl)>0$ and if one of the following holds: either\\
$(a)$ there exists a positive solution $(\ov{\al},\ov{\psi})$ at $\lm=\ov{\lm}$, or\\
$(b)$ if $(\all,\pl)$ is a solution with $\lm\leq\ov{\lm}$, then $\all>0$,\\
then {\rm $\prl$} has a unique solution for any $\lm\leq \ov{\lm}$. In particular the set of solutions in $[0,\ov{\lm}]$
is a real analytic simple curve of positive solutions $\lm\to (\all,\pl)$, $\lm\in [0,\ov{\lm}\,]$ and
$\all\geq \al_{\sscp \ov{\lm}}$, for any $\lm\leq\ov{\lm}$.
\ele
\proof By Lemma 3.4 in \cite{BeBr}, for any $\ov{I}>0$ we already have a uniform estimate for the $L^{\ii}(\om)$ norm of any
solution of \fbi\, but since the change of variables from $(\gal,\val)$ is singular at $I=0$ (see Remark \ref{remeq}),
we cannot use those estimates directly for $\prl$.
\ble\label{lemE1} Let $p\in[1,p_{N})$. For any $\ov{\lm}>0$ there
exists a positive constant $C_1=C_1(r,\om,\ov{\lm},p,N)$ depending only on $\om$, $\ov{\lm}$, $p$, $N$
and $r\in [0,1)$ such that $\|\pl\|_{C^{2,r}_0(\ov{\om})}\leq C_1$ for any solution
$(\all,\pl)$ of {\rm $\prl$} with $\lm\in [0,\ov{\lm}\,]$.
\ele
\proof
Since $\ino (\all+\lm\pl)^p=1$ it is well known (\cite{St4}) that for any $t\in [1,\frac{N}{N-1})$ there exists $C=C(t,N,\om)$ such that $\|\pl\|_{W_0^{1,t}(\om)}\leq C(t,N,\om)$ for any solution of $\prl$.
Thus, by the Sobolev inequalities, for any $1\leq s< \frac{N}{N-2}$ we have $\|\pl\|_{L^s(\om)}\leq C(s,N,\om)$,
for some $C(s,N,\om)$. As a consequence, since $\all\leq 1$ and $p<p_{N}$, then there exists some $m>1$ depending on $p$ and $N$ such that
$\|(\all+\lm\pl)\|_{L^m(\om)}\leq C(p,N,\ov{\lm},s,\om)$,
for any $\lm\leq \ov{\lm}$,
for some $C(p,N,s,\ov{\lm},\om)$ and then, by standard elliptic estimates,
we conclude that $\|\pl\|_{W_0^{2,m}(\om)}\leq C_0(p,N,s,\ov{\lm},\om)$,
for any $\lm\leq \ov{\lm}$, for some
$C_0(p,N,s,\ov{\lm},\om)$.
At this point, since $\om$ is of class $C^3$, the conclusion follows by standard elliptic estimates
and a bootstrap argument.\finedim

\bigskip

{\it Proof of $(i)$.} We argue by contradiction and assume without loss of generality
that a positive solution $(\ov{\al},\ov{\psi})$ exists for $\lm=\ov{\lm}$. Since $\sg_1(\all,\pl)>0$ for any
$\lm\in [0,\ov{\lm}\,]$ by assumption, then by Lemma \ref{lem1.1} we deduce that there exists a small interval around $\ov{\lm}$
where the set of solutions of $\prl$ is a real analytic curve of positive solutions. By Proposition \ref{pr3.2.best} we deduce that
$\all$ is strictly decreasing and therefore in particular that $\all\geq \ov{\al}>0$ in the given left neighborhood.
However we have $\sg_1(\all,\pl)>0$ for any $\lm\in [0,\ov{\lm}\,]$ and by Lemma \ref{lemE1},
Lemma \ref{lem1.1} and Proposition \ref{pr3.2.best}
we can continue this small curve of positive solutions without bifurcation points backward to the left as a curve of positive
solutions defined in $\lm\in (-\dt,\ov{\lm}\,]$ for some small $\dt>0$.
At this point, if at any point in $(0,\ov{\lm}\,]$ a positive solution would exist
other than those on this curve, then it would generate in the same way another curve
of positive solutions in $(-\dt,\ov{\lm}\,]$.
Obviously for both curves at $\lm=0$ we have $(\all,\pl)=(\al_{\sscp 0},\psi_{\sscp 0})$ which is the
unique solution of $\prl$ for $\lm=0$. This is obviously impossible since then $(\al_{\sscp 0},\psi_{\sscp 0})$ would be a bifurcation
point, in contradiction with Lemma \ref{lem1.1}.\\

{\it Proof of $(ii)$.}
If $(a)$ holds we just know that there exists a positive
solution $(\ov{\al},\ov{\psi})$ at $\lm=\ov{\lm}$. Therefore, the same argument adopted in the first part shows that there
is a real analytic curve of positive solutions emanating from $(\ov{\al},\ov{\psi})$ and defined in $\lm\in (-\dt,\ov{\lm}\,]$ for some small $\dt>0$.
In particular, since $\sg_1(\all,\pl)>0$ for any $\lm \in [0,\ov{\lm}\,]$ then by Proposition \ref{pr3.2.best} we deduce that
$\all\geq \ov{\al}$ along this curve.
By the first part we have that there exists at most one positive solution in $[0,\ov{\lm}]$ and we are done in this case.\\
If $(b)$ holds we argue the other way around.
By Lemma \ref{lem1.1} we known that there exists $0<\eps<\ov{\lm}$ small enough such that $\prl$ has at least one solution for any
$\lm\in [0,\eps]$. Since by assumption
$\sg_1(\all,\pl)>0$ for any $\lm\in [0,\ov{\lm}\,]$ and since $(b)$ holds, then by Lemma \ref{lemE1} and Lemma \ref{lem1.1}
we can continue this curve forward to the right as a real analytic curve of positive solutions defined in $\lm\in [0,\ov{\lm}\,]$.
In particular by Proposition \ref{pr3.2.best} we deduce that
$\all\geq \al_{\ov{\lm}}$ along this curve.
By the first part we have that there exists at most one positive solution in $[0,\ov{\lm}]$ and we are done in this case as well.
\finedim

\bigskip
\bigskip

By using Lemma \ref{thm05}, the proof of Theorem \ref{thmLE} is reduced to some uniform
a priori estimates for $\sg_{1}(\all,\pl)$ and $\all$, which we split in various Propositions.

\bpr\label{preigen.LE}
If $(\all,\xil)$ is a solution of {\rm \rife{ul.1}} {\rm (}or equivalently if $(\all,\pl)$ is a solution of {\rm $\prl$}{\rm )} and
$$
\lm p\leq  \Lambda(\om,2p),
$$
then,
{\rm
$$
\nu_{1,\lm}=\inf\limits_{w\in H^{1}_0(\om)\setminus\{0\}}
\dfrac{\ino |\nabla w|^2-\lm p \ino \rlq w^2  }{\ino \rlq w^2}\geq 0
$$
}
\epr
\proof
Obviously we can assume that $\lm>0$. Let $w \in C^{1}_0(\ov{\om}\,)$, $w\equiv \!\!\!\!\!/ \;0$ then we have,
$$
\dfrac{\ino |\nabla w|^2}{\ino \rlq w ^2}\geq \dfrac{1}{ \left(\ino (\all+{\xil})^p\right)^{\frac1q}}
\dfrac{\ino |\nabla w|^2}{\left(\ino w^{2p}\right)^{\frac1p}}=
\dfrac{\ino |\nabla w|^2}{\left(\ino w^{2p}\right)^{\frac1p}}\geq {\Lambda(\om,2p)}
$$
which immediately implies that,
$$
\nu_{1,\lm}\geq {\Lambda(\om,2p)}-\lm p\geq 0,\;\forall\,\lm p\leq\Lambda(\om,2p).
$$
\finedim

\bigskip

\bpr\label{lemsob}
For fixed $p\in[1,p_N)$, it holds $\lm^*(\om,p)\geq  \frac1p \Lambda(\om,2p)$ and the equality holds if and only if $p=1$.
\epr
\proof
It is well known that for $p=1$ (\cite{Te2}) it holds
$\lm^*(\om,1)=\lm^{(1)}(\om)$, and since  $\lm^{(1)}(\om)\equiv \Lambda(\om,2)$, then the equality holds for $p=1$.
Therefore, we will just prove that if $p>1$ then $\lm^*(\om,p)$ is well defined and satisfies
$\lm^*(\om,p)>\frac1p \Lambda(\om,2p)$. To simplify the exposition let us denote $\lm^*=\lm^*(\om,p)$.
By Lemma \ref{lmsmall} in Appendix \ref{appF}, $\lm^*$ is well defined.
We argue by contradiction and suppose that $\lm^*\leq \frac1p\Lambda(\om,2p)$.
It follows from Proposition \ref{preigen.LE}, \rife{xi1s}
and Lemma \ref{thm05} that for any $\lm<\lm^*$ there exists a unique solution of $\prl$ and that
these solutions form a real analytic simple curve of positive solutions. By Lemma \ref{lemE1}
we can pass to the limit and obtain at least one solution $(\al_*,\psi_*)$ of $\prl$ for $\lm=\lm^*$.
If $\al_*=0$ then
$u_*=\lm^*\psi_*$ would be a solution of
\beq\label{us}
\graf{-\Delta u_* = \lm^* (u_*)^p \quad \mbox{in} \;\;\om\\ \\
\bigintss\limits_{\om}(u_*)^p=1\\ \\
u_*>0 \;\;\mbox{in}\;\;\om, \quad u_*=0 \;\; \mbox{on}\;\;\pa\om.
}
\eeq
Now the linearization of \rife{us} where one just disregards the integral constraint takes the form
\beq\label{phis}
\graf{-\Delta \phi = \lm^* p(u_*)^{p-1}\phi  \quad \mbox{in} \;\;\om\\ \\
\phi=0 \quad \mbox{on}\;\;\pa\om,
}
\eeq
and it is readily seen that $u_*$ is a positive strict subsolution of \rife{phis} and then we deduce by standard arguments that the first
eigenvalue of \rife{phis} is negative. In particular we infer that
$$
0>
\inf\limits_{w\in H^{1}_0(\om)\setminus\{0\}}
\dfrac{\ino |\nabla w|^2}{\ino (u_*)^{p-1} w^2}-\lm^* p=
$$
$$
\inf\limits_{w\in H^{1}_0(\om)\setminus\{0\}}
\dfrac{\ino |\nabla w|^2-\lm^* p \ino (\rh_{\sscp\lm^*})^{\frac1q} w^2  }{\ino (\rh_{\sscp\lm^*})^{\frac1q} w^2}=
\nu_{1,\lm^*}(u_*),
$$
which contradicts Proposition \ref{preigen.LE}. Therefore, we must have
$\al_*>0$ and in particular, by the monotonicity of $\all$ (Proposition \ref{pr3.2.best}),
we also have $\inf\limits_{[0,\lm^*)}\all>0$. Using
Proposition \ref{preigen.LE}, \rife{xi1s} and Lemma \ref{lem1.1} we see that there exists $\epsilon$ small enough
such that we can continue the simple curve of unique positive solutions defined on $[0,\lm^*]$ to a larger curve
$\mathcal{G}_{\epsilon}$ defined in
$[0,\lm^*+\epsilon)$ such that in particular, by continuity, $\inf\limits_{(\all,\pl)\in \mathcal{G}_{\epsilon}}\all>0$.
Therefore, by the definition of $\lm^*$, either there exists another solution of $\prl$ for $\lm=\lm^*$ with $\all=0$
or there exists a sequence $(\al_n,\psi_n)$ of solutions of $\prl$ for
$\lm_n\to (\lm^*)^+$ such that $\al_n \to 0$. In the latter case, by Lemma \ref{lemE1}
we can pass to the limit and obtain at least one solution $(\al^*,\psi^*)$ of $\prl$ for $\lm=\lm^*$ with $\al^*=0$.
In particular we deduce that
there exists at least another solution $u^*$, distinct from $u_*$, of \rife{us}.
At this point the same argument adopted above shows that
$\nu_{1,\lm^*}(u^*)<0$ which contradicts once more Proposition \ref{preigen.LE}.
\finedim

\bigskip
\bigskip

At this point we are ready to prove Theorem \ref{thmLE}.\\
{\em Proof of Theorem \ref{thmLE}}.\
By Proposition \ref{lemsob} we have $\lm^*(\om,p)\geq \frac1p\Lambda(\om,2p)$ where the equality holds if and only if $p=1$.
Whence $\all>0$ if $\lm\in [0,\frac1p\Lambda(\om,2p))$.
Therefore,
in view of \rife{xi1s} and Proposition \ref{preigen.LE}, in particular
we see that for any $\lm\in [0,\frac1p\Lambda(\om,2p))$ it holds
$\sg_{1}(\all,\pl)>0$ and $\all>0$ for any solution of $\prl$.
As a consequence by Lemma \ref{thm05} and by Propositions \ref{pr-enrg} and \ref{pr3.2.best},
we deduce  that for any $\lm\in [0,\frac1p\Lambda(\om,2p))$
there exists a unique solution of $\prl$ and that these solutions form a simple real analytic curve of positive
solutions along which we have,
$$
\frac{d \all}{d\lm}<0,\quad \frac{d\el}{d\lm}>0,
\quad \forall\,\lm\in [0,\frac1p\Lambda(\om,2p)).
$$

Next, by well known rearrangement estimates (\cite{Tal2}) we have,
\begin{equation} \label{e0}
E_{0}(\om)=\frac12\ino \ino G_{\om}(x,y)\,dxdy\leq
\frac12\int_{\mathbb{D}_{\sscp N}} \int_{\mathbb{D}_{\sscp N}} G_{\mathbb{D}_{\sscp N}}(x,y)\,dxdy=
{\textstyle\dfrac{|B_1|^{\scp -\frac 2 N}}{4(N+2)}}.
\end{equation}
The fact that $\all=1+\mbox{\rm O}(\lm),\;\; \pl=\psi_{\sscp 0}+\mbox{\rm O}(\lm),\mbox{ as }\lm\to 0^+$ was already part
of Lemma \ref{lem1.1} and then, since $(\all,\pl)$ is a real analytic curve, it is easy to prove that,
$\el=E_{0}(\om)+\mbox{\rm O}(\lm),\mbox{ as }\lm\to 0^+$.\\
By Lemma \ref{lemE1} we can pass to the limit as $\lm\to (\frac1p \Lambda(\om,2p))^{-}$ along $\mathcal{G}(\om)$ to obtain a
solution $(\ov{\al},\ov{\psi})$ of $\prl$ for $\lm=\frac1p \Lambda(\om,2p)$ and
we deduce that $\mathcal{G}(\om)$ can be extended by continuity at $\lm=\frac1p \Lambda(\om,2p)$.
If $p=1$ then $\lm^*(\om,1)=\lm^{(1)}(\om)\equiv \Lambda(\om,2)$ and $\ov{\al}=0$ (\cite{Te2}). Otherwise, since
$\lm^*(\om,p)>\frac1p \Lambda(\om,2p)$ for $p>1$, we must have $\ov{\al}>0$.\finedim

\bigskip
\bigskip

The proof of Theorem \ref{thmvar} is based on the equivalence of the variational
formulations of $\prl$ and \fbi, see Appendix \ref{appD}. Here
\beq\label{lmsdef}
\lm^{**}(\om,1)=\lm^{(1)}(\om) \quad \mbox{and} \quad \lm^{**}(\om,p)=\left(I^{**}(\om,p)\right)^{\frac1q}, \ p>1,
\eeq
denote the threshold for variational solutions of $\prl$, see Theorem A in the introduction and
Corollary \ref{eqlms} in Appendix \ref{appD}.\\

{\it Proof of Theorem \ref{thmvar}.}
By Theorem A in the introduction and Corollary \ref{eqlms} in Appendix \ref{appD}
we know that for variational solutions $\all>0$ if and only if
$\lm<\lm^{**}(\om,p)$.
Since $\lm^{*}(\om,p)$ is by definition the threshold value for any solution, it cannot be larger than $\lm^{**}(\om,p)$,
which proves $\lm^{**}(\om,p)\geq\lm^{*}(\om,p)$ in \rife{istest}. The fact that for $|\om|=1$ it holds
$\Lambda(\om,p)\geq \Lambda(\mathbb{D}_{\sscp 2},p)$ is well known (\cite{CRa}) which together with the inequality
$\lm^{*}(\om,p)\geq \frac1p\Lambda(\om,p)$ of Theorem \ref{thmLE} concludes the proof of \rife{istest}.\\
As a consequence of \rife{istest} we deduce that
any variational solution with $\lm< \frac1p\Lambda(\om,p)$ is a positive solution. Therefore,
it follows from Theorem \ref{thmLE} that it is the unique positive solution in this range.
In particular the set of variational solutions coincides with the set of positive solutions in this range.
\finedim

\bigskip
\bigskip

\section{\bf The case of the two-dimensional ball}\label{sec6}
We are concerned here with the case $\om=\mathbb{D}_{\sscp 2}$ and $p\in[1,+\infty)$. By \cite{BSp} we know that the solutions are unique and we denote with $\mathcal{G}^{*}(\mathbb{D}_{\sscp 2})$ the set of unique solutions of
$\prl$ on $\mathbb{D}_{\sscp 2}$. Recall also \rife{lmsdef} and
Corollary \ref{eqlms} in Appendix \ref{appD}.

\medskip

{\it Proof of Theorem \ref{tdisk}.} For $p=1$ we have by Theorem \ref{thmLE} that
${\lm}^*(\mathbb{D}_{\sscp 2},1)=\lm^{(1)}(\mathbb{D}_{\sscp 2})\equiv \Lambda(\mathbb{D}_{\sscp 2},2)$ and
 the set of solutions $\mathcal{G}^*(\mathbb{D}_{\sscp 2})$ is a real analytic curve in
$[0,{\lm}^*(\mathbb{D}_{\sscp 2},1))$
which is continuous in $[0,{\lm}^*(\mathbb{D}_{\sscp 2},1)]$ and coincides with $\mathcal{G}(\mathbb{D}_{\sscp 2})$. Here we recall that  $\mathcal{G}(\mathbb{D}_{\sscp 2})$ is the set of solutions in Theorem \ref{thmLE}.
In particular the unique (\cite{BeBr,Pudam,Te2}) solution $(\ov{\al},\ov{\psi})$ of $\prl$ with
$\lm={\lm}^*(\mathbb{D}_{\sscp 2},1)$ satisfies $\ov{\al}=0$ and
the monotonicity properties of $\all$ and $\el$ of Theorem \ref{thmLE} hold in $[0,{\lm}^*(\mathbb{D}_{\sscp 2},1))$.\\

Therefore, we are left with the case $p>1$ where we recall that, from Theorem \ref{thmLE},
we have ${\lm}^{*}(\mathbb{D}_{\sscp 2},p)>\frac1p\Lambda(\mathbb{D}_{\sscp 2},2p)$. We will use the fact that, by the
uniqueness of solutions (\cite{BSp}), we have
${\lm}^{*}(\mathbb{D}_{\sscp 2},p)={\lm}^{**}(\mathbb{D}_{\sscp 2},p)$ (see Corollary \ref{eqlms}) and
to simplify the notations in the rest of this proof we set $\lm^*={\lm}^{*}(\mathbb{D}_{\sscp 2},p)$.\\
By Proposition \ref{preigen.LE}, by \rife{xi1s} and Lemma \ref{lem1.1} and since
${\lm}^{*}>\frac1p\Lambda(\mathbb{D}_{\sscp 2},2p)$,
we can continue the curve $\mathcal{G}(\mathbb{D}_{\sscp 2})$ in a right neighborhood of
$\frac1p\Lambda(\mathbb{D}_{\sscp 2},2p)$ as a
larger simple curve of positive solutions with no bifurcation points $\mathcal{G}_{\mu}$, such that
$\mathcal{G}(\mathbb{D}_{\sscp 2})\subset \mathcal{G}_{\mu}$, $\mu>\frac1p\Lambda(\mathbb{D}_{\sscp 2},2p)$, and by continuity
$\all>0$ and $\sg_1(\all,\pl)>0$ for any $\lm<\mu$.
Therefore, it is well defined,
$$
{\lm}_1={\lm}_1(\mathbb{D}_{\sscp 2}):=\sup\left\{\mu>\frac1p\Lambda(\mathbb{D}_{\sscp 2},2p)\,:\,\all>0\mbox{ and }\sg_1(\all,\pl)>0,\;
\forall\,(\all,\pl)\in \mathcal{G}_{\mu},\,\forall \lm<\mu\right\}.
$$
There are three possibilities: either\\
- $\lm_1(\mathbb{D}_{\sscp 2})=+\ii$, or\\
- $\lm_1(\mathbb{D}_{\sscp 2})<+\ii$ and either
$\inf\limits_{\lm\in[0,\lm_1(\mathbb{D}_{\sscp 2}))}\all=0$ or $\inf\limits_{\lm\in[0,\lm_1(\mathbb{D}_{\sscp 2}))}\all>0$.\\
Since $\lm^{**}(\mathbb{D}_{\sscp 2},p)<+\ii$ (see Appendix
\ref{appD}), then $\lm_1(\mathbb{D}_{\sscp 2})<+\ii$, which rules out the first possibility.
By Propositions \ref{pr-enrg} and \ref{pr3.2.best}
$\frac{d\all}{d\lm}<0$ and $\frac{d\el}{d\lm}>0$ continue to hold whenever $\all>0$
and $\sg_1(\all,\pl)>0$ are satisfied, whence for any $\lm<\lm_1(\mathbb{D}_{\sscp 2})$ as well. Now if $\lm_1(\mathbb{D}_{\sscp 2})<+\ii$ and
$\inf\limits_{\lm\in[0,\lm_1(\mathbb{D}_{\sscp 2}))}\all=0$ we have
by definition that $\all>0$ and $\sg_1(\all,\pl)>0$ for any $\lm<\lm_1(\mathbb{D}_{\sscp 2})$, and then by Lemma \ref{lemE1},
the uniqueness
of solutions on $\mathbb{D}_{\sscp 2}$ (\cite{BSp}) and the monotonicity of $\all$ (Proposition \ref{pr3.2.best})
it is not difficult to see that we can continue the curve up to $\lm_1(\mathbb{D}_{\sscp 2})$ and that
$\al_{\lm_1(\mathbb{D}_{\sscp 2})}=0$ holds. In particular
$\lm_1(\mathbb{D}_{\sscp 2})=\lm^*$ and in this case we are done.\\
We are left with the discussion of the case,
$\lm_1(\mathbb{D}_{\sscp 2})<+\ii$ and $\inf\limits_{\lm\in[0,\lm_1(\mathbb{D}_{\sscp 2}))}\all>0$.\\
First of all, by definition and by
Lemmas \ref{lem1.1} and \ref{lemE1}, it is not difficult to see that in this case we have
$\sg_1(\al_1,\psi_1)=0$,  $\al_1=\al_{\sscp \lm_1(\mathbb{D}_{\sscp 2})}>0$ and $\psi_1=\psi_{\sscp \lm_1(\mathbb{D}_{\sscp 2})}$.
By Theorem \ref{spectral1} we see that on $\mathbb{D}_{\sscp 2}$ the transversality condition, needed to apply Proposition \ref{pr3.1},
is always satisfied. Therefore, it follows from
Proposition \ref{pr3.1} that we can continue $\mathcal{G}_{\lm_1(\mathbb{D}_{\sscp 2})}$ to a real analytic
parametrization without bifurcation points,
$$
\mathcal{G}^{(s_1+\dt_1)}=\left\{[0,s_1+\dt_1)\ni s\mapsto (\lm(s),\al(s),\psi(s))\right\},
$$
where, for some $s_1>0$ and $\dt_1>0$, we have that for any
$s\in [0,s_1+\dt_1)$, $(\al(s),\psi(s))$ is a positive solution with $\lm=\lm(s)$ and $\lm(s)=s$ for $s\leq s_1$. We claim
that $\lm(s)$ is monotonic increasing along the branch. Indeed, we recall from \rife{2907.5} in
Proposition \ref{pr3.1} that
$$
<[\phi_1]_{\sscp \lm_1(\mathbb{D}_{\sscp 2})},\psi_1>_{\sscp \lm_1(\mathbb{D}_{\sscp 2})}\neq 0
\mbox{ and $<[\phi_1]_{\sscp \lm_1(\mathbb{D}_{\sscp 2})},\psi_1>_{\sscp \lm_1(\mathbb{D}_{\sscp 2})}$  has the same sign as }
<\phi_1>_{\sscp \lm_1(\mathbb{D}_{\sscp 2})}.
$$
Therefore, there is no loss of generality in assuming
$<[\phi_1]_{\sscp \lm_1(\mathbb{D}_{\sscp 2})},\psi_1>_{\sscp \lm_1(\mathbb{D}_{\sscp 2})}>0$.
At this point we use \rife{2907.10}
in Proposition \ref{pr3.1}, that is, putting $\sg_1(s)=\sg_1(\al(s),\psi(s))$,
$$
\dfrac{\sg_1(s)}{\lm^{'}(s)}=\dfrac{p<[\phi_1]_{\sscp \lm_1(\mathbb{D}_{\sscp 2})},\psi_1>_{\sscp \lm_1(\mathbb{D}_{\sscp 2})}+\mbox{\rm o}(1)}
{<[\phi_1]_{\sscp \lm_1(\mathbb{D}_{\sscp 2})}^2>_{\sscp \lm_1(\mathbb{D}_{\sscp 2})}+\mbox{\rm o}(1)},\mbox{ as }s\to 0,
$$
and infer that $\sg_1(s)\neq 0$ and has the same sign of $\lm^{'}(s)$ for $s$ small enough in a deleted neighborhood of $s=0$.
Indeed the eigenvalues $\sg_k(s)=\sg_k(\al(s), \psi(s))$
of $L_{\ssl(s)}$ are real analytic functions of $s$ (\cite{but}) and then the level sets of a fixed $\sg_k(s)$
cannot have accumulation points unless $\sg_k(s)$ is constant on $(0,s_1+\dt_1)$.
However we can rule out the latter case since for $\lm\leq \frac1p\Lambda(\om,2p)$ we always have $\sg_1(\all,\pl)>0$ and then
no $\sg_n(s)$ can vanish identically.\\
At this point we use again that, by the
uniqueness of solutions (Theorem A in the introduction) any solution of $\prl$ is also a variational solution
(see Appendix \ref{appD}).
Therefore, $(\al(s),\psi(s))$ is a positive variational solution and we infer from Lemma \ref{prmin} that $\sg_1(s)\geq 0$ for
$s-s_1>0$ small enough. Since $\sg_1(s)\neq 0$ for $s\neq 0$
small enough we deduce that $\sg_1(s)>0$ for $s-s_1>0$ small and then in particular that $\lm^{'}(s)>0$. Therefore,
$\lm^{'}(s)>0$ for $s$ in a small enough right neighborhood of $s_1$ and the
curve $(\all,\pl)$ bends right. In particular this implies that, with
the notations of Proposition \ref{pr-enrg}, we can evaluate the following limit,
$$
\lim\limits_{\lm\to \lm_1(\mathbb{D}_{\sscp 2})}\frac{d \el}{d\lm}=
\lim\limits_{\lm\to \lm_1(\mathbb{D}_{\sscp 2})}p\,\ml \left(<[\pl]^2_{\ssl}>_{\ssl}+\lm<[\pl]_{\ssl},\eta_{\ssl}>_{\ssl}\right)=
$$
$$
\mbox{O}(1)+p\,m_{\lm_1(\mathbb{D}_{\sscp 2})}\lim\limits_{\lm\to \lm_1(\mathbb{D}_{\sscp 2})}\lm<[\pl]_{\ssl},\eta_{\ssl}>_{\ssl}=
\mbox{O}(1)+p\,m_{\lm_1(\mathbb{D}_{\sscp 2})} \lim\limits_{s\to 0}\frac{1}{\lm^{'}(s)}\lm(s)<[\psi(s)]_{\sscp \lm(s)},\psi^{'}(s)>_{\sscp \lm(s)}=
$$
$$
\mbox{O}(1)+p\,m_{\lm_1(\mathbb{D}_{\sscp 2})}\lm_1(\mathbb{D}_{\sscp 2})
\lim\limits_{s\to 0}\frac{1}{\lm^{'}(s)}
\left(<[\psi_1]_{\sscp \lm_1(\mathbb{D}_{\sscp 2})},\psi^{'}(0)>_{\sscp \lm_1(\mathbb{D}_{\sscp 2})}+\mbox{o}(1)\right)=
$$
\beq\label{flexen}
\mbox{O}(1)+p\,m_{\lm_1(\mathbb{D}_{\sscp 2})}\lm_1(\mathbb{D}_{\sscp 2})
\lim\limits_{s\to 0}\frac{1}{\lm^{'}(s)}
\left(<[\psi_1]_{\sscp \lm_1(\mathbb{D}_{\sscp 2})},\phi_1>_{\sscp \lm_1(\mathbb{D}_{\sscp 2})}+\mbox{o}(1)\right)=+\ii.
\eeq

Similarly, with the notations of Proposition \ref{pr3.2.best}, we have
$$
\lim\limits_{\lm\to \lm_1(\mathbb{D}_{\sscp 2})}\frac{d \all}{d\lm}=
-\lim\limits_{\lm\to \lm_1(\mathbb{D}_{\sscp 2})}\left(<\pl>_{\ssl}+\lm<\eta_{\ssl}>_{\ssl}\right)=
$$
$$
\mbox{O}(1)-\lim\limits_{\lm\to \lm_1(\mathbb{D}_{\sscp 2})}\lm<\eta_{\ssl}>_{\ssl}=
\mbox{O}(1)-\lm_1(\mathbb{D}_{\sscp 2})\lim\limits_{s\to 0}\frac{1}{\lm^{'}(s)}<\psi^{'}(s)>_{\sscp \lm(s)}=
$$
\beq\label{flexal}
\mbox{O}(1)-\lm_1(\mathbb{D}_{\sscp 2})\lim\limits_{s\to 0}\frac{1}{\lm^{'}(s)}\left(<\phi_1>_{\lm_1(\mathbb{D}_{\sscp 2})}+\mbox{o}(1)\right)=-\ii.
\eeq

Summarizing we conclude that there exists $\ov{\lm}_1$ such that we can continue $\mathcal{G}_{\lm_1(\mathbb{D}_{\sscp 2})}$ to a larger
branch $\mathcal{G}_{\ov{\lm}_1}$ such that:\\
$(A1)_1$ for any $\lm \in [0,\ov{\lm}_1)$, $(\all,\pl)$ is a positive solution;\\
$(A1)_2$ $\ov{\lm}_1>\lm_1(\mathbb{D}_{\sscp 2})>\frac1p\Lambda(\mathbb{D}_{\sscp 2},2p)$;\\
$(A1)_3$ the inclusion $\left\{(\all,\pl),\,\lm \in [0,\lm_1(\mathbb{D}_{\sscp 2})] \right\}\subset \mathcal{G}_{\ov{\lm}_1}$,
holds,\\
$(A1)_4$ $\all>0,\, \forall\,\lm \in [0,\ov{\lm}_1]$,\\
$(A1)_5$ $\sg_1(\all,\pl)>0,\, \forall\,\lm\in [0,\ov{\lm}_1]\setminus\{\lm_1(\mathbb{D}_{\sscp 2})\}$.\\
$(A1)_6$ $\frac{d \all}{d\lm}<0$ and $\frac{d\el}{d\lm}>0$ for $\lm \in [0,\ov{\lm}_1]\setminus\{\lm_1(\mathbb{D}_{\sscp 2})\}$
and $\frac{d \all}{d\lm}\to -\ii$, $\frac{d \el}{d\lm}\to +\ii$ as $\lm\to \lm_1(\mathbb{D}_{\sscp 2})$.\\
Clearly, $(A1)_6$ follows from $(A1)_4$, $(A1)_5$, Propositions \ref{pr-enrg}, \ref{pr3.2.best} and \rife{flexen}, \rife{flexal}.\\

Therefore, we can argue by induction and for $k\geq 2$ define,
$$
{\lm}_k(\mathbb{D}_{\sscp 2}):=\sup\left\{\mu>\lm_{k-1}(\mathbb{D}_{\sscp 2})\,:\,\all>0\mbox{ and }\sg_1(\all,\pl)>0,\;
\forall\,(\all,\pl)\in \mathcal{G}_{\mu},\,\forall \lm_{k-1}(\mathbb{D}_{\sscp 2})<\lm<\mu\right\},
$$
where as above we know that $\lm_k(\mathbb{D}_{\sscp 2})<+\ii$. If there exists some $k\geq 2$
such that we also have $\inf\limits_{s\in(0,s_k)}\al(s)=0$, then by arguing as above we are done.
Otherwise, the above procedure yields
a sequence $\ov{\lm}_k$ such that, for any $k\in \mathbb{N}$ we have,
$\ov{\lm}_k>\ov{\lm}_{k-1}>\cdots>\ov{\lm}_{2}>\ov{\lm}_{1}$ and
we can continue $\mathcal{G}_{\lm_1(\mathbb{D}_{\sscp 2})}$ to a larger branch $\mathcal{G}^*(\mathbb{D}_{\sscp 2})$, such that, for any $k\in \N$,
it holds:\\
$(Ak)_1$ for any $\lm \in [0,\ov{\lm}_k)$, $(\all,\pl)$ is a positive solution;\\
$(Ak)_2$ $\ov{\lm}_k>\lm_k(\mathbb{D}_{\sscp 2})>
{\lm}_{k-1}(\mathbb{D}_{\sscp 2})>\cdots>{\lm}_{2}(\mathbb{D}_{\sscp 2})>{\lm}_{1}(\mathbb{D}_{\sscp 2})>\frac1p\Lambda(\mathbb{D}_{\sscp 2},2p)$;\\
$(Ak)_3$ the inclusion $\left\{(\all,\pl),\,\lm \in [0,\lm_k(\mathbb{D}_{\sscp 2})] \right\}\subset \mathcal{G}^{*}(\mathbb{D}_{\sscp 2})$,
holds;\\
$(Ak)_4$ $\all>0,\, \forall\,\lm \in [0,\ov{\lm}_k]$;\\
$(Ak)_5$ $\sg_1(\all,\pl)>0,\, \forall\,\lm\in [0,\ov{\lm}_k]\setminus\{\lm_k(\mathbb{D}_{\sscp 2})\}_{k\in\N}$;\\
$(Ak)_6$ $\frac{d \all}{d\lm}<0$ and $\frac{d\el}{d\lm}>0$ for $\lm \in [0,\ov{\lm}_1]\setminus\{\lm_k(\mathbb{D}_{\sscp 2})\}_{k\in\N}$
and $\frac{d \all}{d\lm}\to -\ii$, $\frac{d \el}{d\lm}\to +\ii$ as $\lm\to \lm_k(\mathbb{D}_{\sscp 2})$.\\

At this point, since $\sup\limits_{k\in \mathbb{N}}{\lm}_k(\mathbb{D}_{\sscp 2})<+\ii$, then
${\lm}_k(\mathbb{D}_{\sscp 2})\to \ov{\lm}$ as $k\to +\ii$. We claim that necessarily $\al_k=\al_{\sscp {\lm}_k(\mathbb{D}_{\sscp 2})}\to 0^+$.
We argue by contradiction and assume that $\al_k\to \ov{\al}\in (0,1]$.
Clearly, by Lemma \ref{lemE1} there exists $\ov{\psi}$ such that
$(\ov{\al},\ov{\psi})$ is a positive solution. Let $\{\ov{\sg}_n\}_{n\in\N}$ be the eigenvalues corresponding to $(\ov{\al},\ov{\psi})$.
First of all we have that $0\in \sigma(L_{\sscp \ov{\lm}})$, where $\sigma(L_{\sscp \ov{\lm}})$ stands for the spectrum of $L_{\sscp \ov{\lm}}$. Indeed, if this was not the case, then, since
$\lm(t_k)\to\ov{\lm}$, by Lemma \ref{lem1.1}, by continuity and since the eigenvalues are isolated, we would have
that there exists a fixed full neighborhood of $0$ with empty intersection with $\sigma(L_{\sscp {\lm_k}})$ for any $k$ large enough.
This is clearly impossible since by construction $0\in \sigma(L_{\sscp {\lm_k}})$ for any $k$.
As a consequence $\ov{\sg}_1=0$. However, since $\ov{\al}\in (0,1]$ by assumption, then
we can apply
Proposition \ref{pr3.1} and in particular conclude that $\ov{\sg}_1(s)$
is a real analytic function of $s$ (\cite{but}) and then its zero level set
cannot have accumulation points unless it vanishes identically. On one side for
$\lm\leq \frac1p\Lambda(\mathbb{D}_{\sscp 2},2p)$ we have $\sg_1(\all,\pl)>0$ and then we deduce
from $(A1)_2$ that no eigenvalues can vanish identically. Therefore, since the eigenvalues are isolated, by continuity
we would have once more that there exists a fixed full neighborhood of $\ov{\sg}_1=0$ with empty intersection with
$\sigma(L_{\sscp {\lm_k}})$ for any $k$ large enough.
This is clearly impossible since by construction $0\in \sigma(L_{\sscp {\lm_k}})$ for any $k$.
Therefore, $\ov{\al}=0$ and $\ov{\lm}=\lm^*$ in this case as well.\\
The energy identity follows by Theorem 1.1 in \cite{BJ2}.
\finedim

\bigskip
\bigskip

\appendix

\section{Variational solutions}\label{appD}
Problem  $\prl$ arises as the Euler-Lagrange equation of the constrained
minimization principle {\bf (VP)} below for the plasma densities {\rm $\rh\in L^{1}(\om)$},
which, for $p>1$, is equivalent to the variational formulation of \fbi.
We shortly discuss here the variational solutions of \fbi\, and $\prl$ and their equivalence and
refer to \cite{BeBr} for a detailed discussion of this point. In this context $\all$ is the Lagrange multiplier related to the
"mass" constraint $\ino \rl=1$ while the Dirichlet energy is the density interaction energy,
$$
\mathcal{E}(\rh)=\frac12\ino \rh G[\rh],
$$
which is easily seen to coincide with $\el$ whenever $\pl=G[\rl]$, that is, $\el=\mathcal{E}(\rl)$.\\
For any
$$
\rh\in\mathcal{P}_{\sscp \om}:=\left\{\rh\in L^{1+\frac1p}(\om)\,|\,\rh\geq 0\;\mbox{a.e. in}\;\om \right\},
$$
and $\lm\geq  0$, we define the free energy,
\beq\label{jeil}
J_{\ssl}(\rh)=
{\scriptstyle \frac{p}{p+1}}\ino (\rh)^{1+\frac{1}{p}}-\frac\lm 2 \ino \rho G[\rho].
\eeq
Let us consider the variational principle,
$$
\mathcal{J}(\lm)=\inf\left\{ J_{\ssl}(\rh)\,:\,\rh\in \mathcal{P}_{\sscp \om}, \ino \rh=1\right\}.\qquad \qquad \mbox{\bf (VP)}
$$
It has been shown in \cite{BeBr, Te2} that for each $\lm>0$ there exists at least one $\rl$ which solves the {\bf (VP)}.
In particular, (\cite{BeBr}) $\al=\all\in\R$ arises as the Lagrange multiplier relative to the constraint $\ino \rl=1$ and if
$\all\geq 0$, then
any minimizer $\rl$ yields a solution $(\all,\pl)$ of $\prl$ where $\pl=G[\rl]$. Any such solution is called a variational solution of $\prl$.\\

Solutions of \fbi\, are also found in \cite{BeBr} as minimizers of the
functional $\Psi_{I}(v)$, $v\in\mathcal{H}_I$ as defined in \rife{psiI}.
For fixed $I>0$, a variational solution of \fbi\, is a solution of \fbi\, which is also a minimizer of
$\Psi_I$ on $\mathcal{H}_I$. It has been shown in \cite{BeBr} that, for $p\in(1,p_N)$, at least one variational solution exists for each $I>0$.
From Theorem A in the introduction there exists $I^{**}(\om,p)>0$ such that if $v$ is a variational solution of
\fbi\, for some $\ga=\gal$, then $\gal>0$ if and only if $I\in (0,I^{**}(\om,p))$.\\
We will use a result in \cite{BeBr} p.421 which shows that, for $p>1$, there is a one to one correspondence between variational
solutions of \fbi\, and minimizers of
$$
\inf\left\{ J_{1}(\widetilde{\rh})\,:\,\widetilde{\rh}\in \mathcal{P}_{\sscp \om}, \ino \widetilde{\rh}=I\right\},
$$
where $J_\lm$ has been defined in \rife{jeil}. We will also use the fact that the minimization of $J_1$
 on $\mathcal{P}_{\sscp \om}\cap\{\ino \widetilde{\rh}=I\}$ is
 equivalent to {\bf (VP)}. This is readily seen from {\bf (VP)} after the scaling $\rh=I^{-1}\widetilde{\rh}$ with
 $\lm=I^{\frac1q}$ and in particular yields the Euler-Lagrange equation,
\beq\label{euler.1}
(\mbox{\rm $\widetilde{\rh}$})^{\frac1p}=(\widetilde{\al} + G[\mbox{\rm $\widetilde{\rh}$}])_+ \mbox{ in }\om.
\eeq
Actually, it is not difficult to see that, if $\widetilde{\al}\geq 0$, then \rife{euler.1} is just the same as $\prl$ with
$\rl=I^{-1}\widetilde{\rh}$, $\lm=I^{\frac1q}$ and $\all=I^{-\frac{1}{p}}\widetilde{\al}$.\\
From \cite{BeBr} p.421, we have: if $\widetilde{v}$ minimizes $\Psi_I$ on $\mathcal{H}_I$,
then $\widetilde{\rh}=(\widetilde{v})^p_+$ is a minimizer of $J_1$
 on $\mathcal{P}_{\sscp \om}\cap\{\ino \widetilde{\rh}=I\}$. Therefore,
 if $(\gal,\widetilde{v}_{\sscp I})$ is a positive variational solution,
 then $\widetilde{\rh}_{\sscp I}=(\widetilde{v}_{\sscp I})^p_+\equiv \widetilde{v}^p_{\sscp I}$ and in particular
$\left.\widetilde{\rh}_{\sscp I}\right|_{\pa \om}=:\gal^p>0$. Therefore, it follows from \rife{euler.1} that
$\widetilde{\al}_{\sscp I}=(\widetilde{\rh}_{\sscp I})^{\frac1p}=\gal>0$, and then scaling back we see that
$(\all,\pl)=(I^{-\frac{1}{p}}\widetilde{\al}_{\sscp I},  G[I^{-1}\widetilde{\rh}_{\sscp I}])$ is a positive variational solution of {\bf (VP)}.\\
On the other side, still from \cite{BeBr} p.421, we have: if $\widetilde{\rh}$ is a minimizer of $J_1$
 on $\mathcal{P}_{\sscp \om}\cap\{\ino \widetilde{\rh}=I\}$, then there exists a unique $\gamma\in \R$ such that
$\widetilde{v}=\gamma+G[\widetilde{\rh}] \in \mathcal{H}_I$ and $\widetilde{v}$ minimizes $\Psi_I$ on $\mathcal{H}_I$.
Now if $(\all,\pl)$ is a positive variational solution of $\prl$ and $\rl=(\all +\lm\pl)^p$ then
$\widetilde{\rh}_{\sscp I}=I\left(\rl\right)_{\lm =I^{\frac1q}}$ minimizes $J_1$ on
$\mathcal{P}_{\sscp \om}\cap\{\ino \widetilde{\rh}=I\}$ and we infer that
$\widetilde{v}_{\sscp I}=\gal+G[\widetilde{\rh}_{\sscp I}]= \widetilde{\rh}^{\frac1p}_{\sscp I}$ is a minimizer
of $\Psi_I$ on $\mathcal{H}_I$ for a unique $\gal$ which satisfies
$\gal=\left. \widetilde{\rh}^{\frac1p}_{\sscp I}\right|_{\pa\om}={I}^{\frac1p}\left(\all\right)_{\lm =I^{\frac1q}}>0$.
Therefore, we conclude that $(\gal,\val)$ is positive variational solution of
\fbi.\\

Summarizing, by using a duality argument introduced in \cite{BeBr}, we have that, for $p>1$, $(\all,\pl)$ is a
positive variational solution of {\rm $\prl$} if and only if $(\gal,\val)$
is a positive variational solution of \fbi. Therefore, as a corollary of Theorem A in the introduction,
we conclude that
\bco\label{eqlms}
Let $p\in[1,p_N)$. Then there exists $\lm^{**}(\om,p)\in (0,+\ii)$ such that $(\all,\pl)$ is a positive variational solution
of {\rm $\prl$} if and only if $\lm\in (0,\lm^{**}(\om,p))$. In particular, for $p>1$, via Remark \ref{remeq} we have
$\lm^{**}(\om,p)=\left(I^{**}(\om,p)\right)^{\frac1q}$. Moreover, $\lm^{*}(\mathbb{D}_{\sscp 2},p)=\lm^{**}(\mathbb{D}_{\sscp 2},p)$.
\eco
The identity $\lm^{*}(\mathbb{D}_{\sscp 2},p)=\lm^{**}(\mathbb{D}_{\sscp 2},p)$ is an immediate consequence of the uniqueness
of solutions on $\mathbb{D}_{\sscp 2}$ (\cite{BSp}), which implies that any solution is variational.

\bigskip
\bigskip

\section{Uniqueness of solutions for $\lm$ small}\label{appF}
The following lemma is proved by a standard application of the contraction mapping principle and we
prove it here just for reader's convenience.
\ble\label{lmsmall}
There exists $\ov{\lm}>0$ such that for any $\lm\in [0,\ov{\lm}\,]$ there exists one and only one solution
$(\all,\pl)$ of {\rm $\prl$}. Moreover, $\all\geq \frac13$ for any $\lm\in [0,\ov{\lm}\,]$.
\ele
\proof The proof is an immediate consequence of the following lemmas. Let us define
$$
B_{\ii}=\left\{u\in L^{\ii}(\om)\,|\,\|u\|_{L^{\ii}(\om)}\leq C_1(r,\om),\, u\geq 0 \mbox{ a.e. in } \om\right\}
$$
where $C_1(r,\om)$ is the constant in Lemma \ref{lemE1} concerning the uniform a priori bounds.

\ble\label{lemE2} Let $p\in[1,p_N)$. There exists $\ov{\lm}>0$ such that for any $\lm\in [0,\ov{\lm}\,]$ and for any $\al\in [0,1]$ there exists
one and only one solution $u=u_{\sscp \lm,\al}\in C^{2,r}_0(\ov{\om}\,)$ of the problem
$$
\graf{-\Delta u =\lm (\al+u)^p\quad \mbox{in}\;\;\om\\ \\
u=0 \quad \mbox{on}\;\;\pa\om\\\\
u\in B_{\ii}.}
$$
Moreover, for fixed $\lm\in [0,\ov{\lm}]$, the map
$[0,1]\ni \al\to u_{\lm}[\al]=u_{\lm,\al}\in B_{\ii}$ is continuous and $u_{0,\al}=0$.
\ele
\proof
For $\ov{\lm}>0$ to be fixed later on and for fixed $\lm\in [0,\ov{\lm}\,]$ and $\al\in [0,1]$ we define
$$
T_{\lm,\al}(u)=\lm G[(\al+u)^p],\quad u\in B_{\ii}.
$$
Clearly
$$
\|T_{\lm,\al}(u)\|_{L^{\ii}(\om)}\leq \lm C_2(p,\om,\ov{\lm}C_1(r,\om)),
$$
whence we have
$T_{\lm,\al}:B_{\ii}\to B_{\ii}$ for any $\lm\leq \frac{1}{C_2}$. Also,
$$
\|T_{\lm,\al}(u)-T_{\lm,\al}(v)\|_{L^{\ii}(\om)}\leq \lm p\|G[(\al+w)^{p-1}|u-v|]\|_{L^{\ii}(\om)}\leq
$$
$$
\lm p\|G[(\al+w)^{p-1}]\|_{L^{\ii}(\om)} \|u-v\|_{L^{\ii}(\om)}\leq \lm p\|G[(1+|w|)^{p-1}]\|_{L^{\ii}(\om)} \|u-v\|_{L^{\ii}(\om)}
\leq
$$
$$
\lm C_3\|u-v\|_{L^{\ii}(\om)}
$$
where $w\in B_{\ii}$ satisfies $u\leq w\leq v$ and $C_3$ depends only by $p,\om$ and $\ov{\lm}C_1(r,\om)$.

Therefore, we also have $\|T_{\lm,\al}(u)-T_{\lm,\al}(v)\|_{L^{\ii}(\om)}\leq \frac12\|u-v\|_{L^{\ii}(\om)}$, for any
$\lm\leq \frac{1}{2C_3}$. As a consequence choosing $\ov{\lm}\leq \min\{\frac{1}{C_2},\frac{1}{2C_3}\}$, then
$T_{\lm,\al}$ is a contraction on $B_{\ii}$.
Whence, in particular, for any fixed $\al\in [0,1]$, we have that for any $\lm\in [0,\ov{\lm}\,]$ there exists a unique solution of $u=T_{\lm,\al}(u)$. The existence and uniqueness claim follows since,
by standard elliptic estimates, $u=u_{\sscp \lm,\al}\in C^{2,r}_0(\ov{\om})$ solves the problem
in the statement of the lemma if and only if $u\in B_{\ii}$ satisfies $u=T_{\lm,\al}(u)$.\\
Concerning the continuity of $u_{\lm}[\al]=u_{\lm,\al}$ for $\al\in [0,1]$, we observe that if $\al_n\to \al$, then
$$
\|u_{\lm}[\al_n]-u_{\lm}[\al]\|_{L^{\ii}(\om)}=\|T_{\lm,\al_n}(u_{\lm,\al_n})-T_{\lm,\al}(u_{\lm,\al})\|_{L^{\ii}(\om)}\leq
$$
$$
\|T_{\lm,\al_n}(u_{\lm,\al_n})-T_{\lm,\al_n}(u_{\lm,\al})\|_{L^{\ii}(\om)}+
\|T_{\lm,\al_n}(u_{\lm,\al})-T_{\lm,\al}(u_{\lm,\al})\|_{L^{\ii}(\om)}\leq
$$
$$
\lm C_3\|u_{\lm}[\al_n]-u_{\lm}[\al]\|_{L^{\ii}(\om)} + p\lm\|G[(s+u_{\lm,\al})^{p-1}]\|_{L^{\ii}(\om)}|\al_n-\al|\leq
$$
$$
\frac{1}{2}\|u_{\lm}[\al_n]-u_{\lm}[\al]\|_{L^{\ii}(\om)}+\ov{\lm}C_3|\al_n-\al|,
$$
which readily implies the claim. Obviously $u_{0,\al}=T_{0,\al}[u_{0,\al}]=0$.
\finedim

\bigskip

For fixed $\lm\in [0,\ov{\lm}\,]$ we consider the continuous map $[0,1]\ni \al\to u_{\lm}[\al]\in B_{\ii}$ where $u_{\lm}[\al]=u_{\lm,\al}$.
Then we have,

\ble\label{lemE3} By taking a smaller $\ov{\lm}$ if necessary, for any fixed $\lm\in [0,\ov{\lm}\,]$ we have:\\
$(i)$ The map $u_{\lm}[\al]$ is monotonic increasing,
$$
u_{\lm}[\al]\leq u_{\lm}[\beta],\quad \forall\,0<\al<\beta\leq 1.
$$
$(ii)$ There exists a unique $\all\in [\frac13,1]$ such that,
$$
\ino (\all+u_{\lm}[\all])^p=1.
$$
\ele
\proof
$(i)$ If $\lm=0$ we have $u_{0,\al}=0$ for any $\al$ and the conclusion is trivial. For any fixed $0<\al<\beta\leq 1$ let $w=u_{\lm}[\beta]-u_{\lm}[\al]\in C^{2,r}_0(\ov{\om}\,)\cap B_{\ii}$, then
$$
-\Delta w=\lm (\beta+u_{\lm}[\beta])^p-\lm(\al+u_{\lm}[\al])^p>
$$
$$
\lm p (\al+u_{\lm}[\al])^{p-1}(u_{\lm}[\beta]-u_{\lm}[\al])=\lm p (\al+u_{\lm}[\al])^{p-1}w,
$$
by the convexity of $f(t)=(\al+t)^p$ for $\al\in (0,1]$. Observe now that the first eigenvalue $\xi_{1,\lm}$
of the linearized operator
$-\Delta w +\lm f^{'}(u_{\lm}[\al])w$, $w \in H^{1}_0(\om)$ satisfies,

\beq\label{xi1n}
\xi_{1,\lm}=\inf\limits_{w\in H^{1}_0(\om)\setminus\{0\}}\dfrac{\ino |\nabla w|^2-\lm \ino f^{'}(u_{\lm}[\al])|w|^2  }{\ino |w|^2}=
\eeq
$$
\inf\limits_{w\in H^{1}_0(\om)\setminus\{0\}}\left(\dfrac{\ino |\nabla w|^2}{\ino |w|^2}-
\dfrac{\lm \ino f^{'}(u_{\lm}[\al])|w|^2  }{\ino |w|^2} \right)\geq
\lm^{(1)}(\om)-\lm p(1+\ov{\lm}C_1(r,\om))^{p-1}\geq \frac{1}{2}\lm^{(1)}(\om),
$$

where possibly we take a smaller $\ov{\lm}$ to guarantee that $p\ov{\lm}(1+\ov{\lm}C_1(r,\om))^{p-1}\leq \frac{1}{2}\lm^{(1)}(\om)$. Here $\lm^{(1)}(\om)$ is the first eigenvalue of $-\Delta$ on $\om$ with Dirichlet boundary conditions and $C_1(r,\om)$ is the constant appearing in Lemma \ref{lemE2}. In particular
 we conclude that $w=u_{\lm}[\beta]-u_{\lm}[\al]\geq 0$, as claimed.\\
 $(ii)$ For $\lm=0$ we have $u_{0,\al}=0$ and then necessarily $\all=1$. For fixed $\lm \in (0,\ov{\lm}\,]$ and by Lemma \ref{lemE2} and $(i)$, the function
$$
g(\al)=\ino (\al+u_{\lm}[\al])^p,\quad \al \in (0,1],
$$
is continuous and strictly increasing. Moreover, possibly
taking $\ov{\lm}$ small enough to guarantee that $(2 C_1\ov{\lm} \,)^p\leq \frac{1}{4}$, we have that,
$$
\limsup\limits_{\al \to 0} g(\al)\leq \limsup\limits_{\al \to 0} \left(2^p{\al^p}+2^p\ino  (u_{\lm}[\al])^p\right)\leq \frac14,
$$
while, for any $\lm\in(0,\ov{\lm}\,]$, we also have,
$$
g(1)=\ino (1+u_{\lm}[\al])^p >1.
$$
Therefore, for any $\lm \in (0,\ov{\lm})$, there exists one and only one $\all$ such that $g(\all)=1$.
On the other side we have,
$$
1 =\ino (\all+u_{\lm}[\all])^p\leq 2^p\all^p +(2 C_1\ov{\lm} \,)^p\leq 2^p\all^p+\frac{1}{4},
$$
whence $\all\geq (\frac{3}{4})^{\frac1p}\frac12>\frac13$.
\finedim

\end{document}